\documentclass[]{siamart190516}
\usepackage{amsmath,amssymb,enumerate,bm,url}
\usepackage{hyperref,cleveref}
\usepackage{microtype}
\usepackage{arydshln}
\usepackage{appendix}
\usepackage{graphicx}
\usepackage{caption}

\usepackage[caption=true]{subfig}
\allowdisplaybreaks

\newcommand{\bs}[1]{\boldsymbol{#1}}

\def \ri {{\rm i}}
\newsiamremark{remark}{Remark}
\newsiamremark{assumption}{Assumption}
\begin{document}

\title{On Spectral Bias Reduction of Multi-scale Neural Networks for Regression Problems \thanks{October 16, 2023, submitted to SIAM J. Numerical Analysis}}
\author{ Bo Wang\thanks{LCSM, Ministry of Education, School of Mathematics and
Statistics, Hunan Normal University, Changsha, Hunan 410081, P. R. China.
The author acknowledges the financial support provided by NSFC (grant
12022104,12371394).}
\and Heng Yuan\thanks{LCSM, Ministry of Education, School of Mathematics and
	Statistics, Hunan Normal University, Changsha, Hunan 410081, P. R. China.}
\and Lizuo Liu\thanks{Department of Mathematics, Southern Methodist
	University, Dallas, TX 75275.}
\and Wenzhong Zhang\thanks{Suzhou Institute for Advanced Research, University of Science and Technology of China, Suzhou, Jiangsu 215000, China. The author acknowledges the financial support provided by NSFC (grant 12201603)}
\and Wei Cai\thanks{Corresponding author, Department of Mathematics, Southern
Methodist University, Dallas, TX 75275 (\texttt{cai@smu.edu}). The author acknowledges the financial support provided by the US National Science
Foundation grant DMS-2207449.}
}

\maketitle

\begin{abstract}
In this paper, we derive diffusion equation models in the spectral domain for the evolution of training errors of two-layer multi-scale deep neural networks (MscaleDNN) \cite{caixu2019,liu2020multi}, designed to reduce the spectral bias of fully connected deep neural networks in approximating oscillatory functions. The diffusion models are obtained from the spectral form of the error equation of the MscaleDNN, derived with a neural tangent kernel approach and gradient descent training and a sine activation function, assuming a vanishing learning rate and infinite network width and domain size. The involved diffusion coefficients are shown to have larger supports if more scales are used in the MscaleDNN, and thus, the proposed diffusion equation models in the frequency domain explain the MscaleDNN's spectral bias reduction capability.  Numerical results of the diffusion models for a two-layer MscaleDNN training match with the error evolution of actual gradient descent training with a reasonably large network width, thus validating the effectiveness of the diffusion models. Meanwhile, the numerical results for MscaleDNN show  error decay over a wide frequency range and confirm the advantage of using the MscaleDNN in approximating functions with a wide range of frequencies.
%
\end{abstract}



\section{Introduction}
Deep learning algorithms  have achieved great success in computer vision \cite{krizhevsky2017imagenet,simonyan2014very,traore2018deep}, natural language processing \cite{young2018recent,otter2020survey,lauriola2022introduction} and many other areas. Their computational power with the help of graphics processing units (GPUs) and capability of handling high dimensional problems have led the computational community to investigate their potentials in applied mathematics research. As a result, a new
research field known as scientific machine learning has became active in the past few
years.

One important task in scientific machine learning is to use deep neural networks (DNNs) to approximate functions or solutions of partial differential equations (PDEs). The idea of using neural networks to solve PDEs  goes back to the 1990's \cite{Dissanayake1994neural,lagaris1998artificial}. In general, four categories of deep PDE solvers have been investigated. The first category is to use deep neural networks to improve classical numerical methods \cite{greenfeld2019learning,hsieh2019learning,um2020solver}. In the second category, the solution operators between infinite-dimensional spaces are approximated by neural networks \cite{anandkumar2020neural,li2020multipole,li2020fourier}. In the third category, the deep neural networks are utilized to approximate
the solutions of PDEs directly such as the physics-informed neural networks (PINNs) \cite{liu2021vpvnet,raissi2018hidden,raissi2019physics}, the deep Ritz method  \cite{yu2018deep,muller2019deep,liao2019deep,hu2022solving}, and Garkerkin methods with weak adversarial networks (WAN) \cite{zang2020weak, Fredlearning}. Lastly, Feynman-Kac formula approaches utilize the connection between linear and nonlinear PDEs and (backward) stochastic differential equations to construct loss functions for the learning algorithms \cite{han2017deep,han2020convergence,han2018solving,beck2019machine,zhang2022fbsde,deepMart}.

Despite their many successes for a wide range of applications, recent studies on the convergence of the deep learning algorithms in theories and practical computations show that standard fully connected DNNs have difficulties in learning high frequency functions, a phenomenon referred as "spectral bias"  \cite{rahaman2019spectral} or  "F-Principle" \cite{xu2020frequency, xu2022overview} in the literature. To overcome this bias, several strategies have been introduced to design neural networks with better frequency resolution, producing promising results.  For solving PDEs, a multi-scale DNN (MscaleDNN) \cite{caixu2019,liu2020multi,wang2020multi,lu2018multimodal,li2020multi,zhang2023corre}, which consists of a series of parallel fully connected sub-neural networks receiving scaled inputs,  has been proposed to learn highly oscillating solutions. Each individual scaled sub-network in the MscaleDNN is designed to approximate a segment of frequency content of the target function, and the effect of the scaling is to convert a specific range of high frequency content to a lower one so that the learning can be accomplished much faster. It was also proposed in \cite{caixu2019, liu2020multi} that the MscaleDNN should use activation functions with a localized frequency profile, such as the sine function and compact supported functions (hat functions and B-splines, etc). In a related work for image and 3D shape reconstruction, the Fourier feature networks, which in fact can be obtained from the MscaleDNN with a sine activation function in their first layer, use the sinusoidal mapping on their inputs and dramatically improve the performance of learning \cite{mildenhall2021nerf,zhong2021cryodrgn,tancik2020fourier}.

To analyze the convergence of deep learning algorithms, the neural tangent kernel (NTK), introduced in \cite{jacot2018neural}, has been a very effective tool to study the evolution of DNNs in function spaces during training \cite{arora2019fine,lee2019wide,ma2020machine,luo2022exact,peng2023non}, and the eigenvector space for the NTK provides many information on the convergence of the DNNs. The convergence and spectral bias of the standard fully connected models and
the improved performance of the afore-mentioned Fourier feature embedded neural networks can be explained by using the NTK. In fact, the NTK theory suggests that standard fully connected DNN have a kernel with a rapid frequency falloff, which prevents them from being able to represent the high-frequency contents of target functions effectively. Fourier feature embedded neural networks have been designed to modify the Fourier spectrum of the NTK so that a faster training convergence for high frequency components can be achieved \cite{mildenhall2021nerf,zhong2021cryodrgn,tancik2020fourier,ParisFourier2020}.

Most of the convergence analysis so far has  been done in the physical domain \cite{tancik2020fourier,luo2022exact,peng2023non,ronen2019convergence, lee2019wide}. Some explicit formulas of NTKs have been  reported for two layers neural networks with the ReLU activation function  \cite{williams1996computing, cho2009kernel, ronen2019convergence,arora2019fine,xie2017diverse,tsuchida2018invariance}. The behaviors of the NTK are usually obtained by analyzing the eigenvalues of the corresponding Gram matrix \cite{tancik2020fourier,peng2023non}.
In this paper, in order to illuminate the mechanism behind the observed reduced spectral bias in the convergence of the MscaleDNN in approximating highly oscillatory functions and PDE solutions \cite{caixu2019,liu2020multi,wang2020multi},
we will derive an error diffusion model, using the NTK approach, in the spectral domain for a two-layered MscaleDNN with a sine activation function  for the case of vanishing learning rate and infinite network width and domain size. Our contribution is three folds: i) we prove that the gradient descent training is equivalent to a diffusion problem in the Fourier spectral domain; ii) the diffusion coefficients can be determined by the Fourier transform of the NTK; iii) the MscaleDNNs with more scales result in diffusion coefficients with larger value and support in the frequency domain. Therefore, our theoretical results provide clear a mathematical explanation why the MscaleDNN can learn much faster over a wider range of frequency. Also, due to the connection between the MscaleDNN and Fourier feature network \cite{tancik2020fourier}, the presented theory can be applied to the latter, as well.

The rest of the paper is organized as follows. In Section 2, a brief review of the MscaleDNN is given. Section 3 derives the diffusion equation models for the training error of high dimensional fitting problem. Analysis of the spectral bias reduction of a two layer MscaleDNN will be done by solving the error diffusion equation models using a Hermite spectral method in Section 4. The numerical results show that the MscaleDNN leads to faster convergence over wider range of frequencies when the number of scales is increased. Finally, section 5 gives a conclusion and  future work.

\section{A review of the multi-scale DNN (MscaleDNN)}
The frequency bias behavior of the deep learning algorithms \cite{rahaman2019spectral, xu2020frequency} has inspired the development and usage of the MscaleDNN in various applications. The MscaleDNN is simply a combination of several fully connected DNNs with different scales on their inputs. It is very convenient to replace a fully connected DNN by a MscaleDNN with equal number of total neurons  in a deep learning algorithm while much better results can be expected. The main idea of the MscaleDNN is to do a radial scaling in the frequency domain such that the learning is performed on functions of scaled-down frequency ranges \cite{liu2020multi,wang2020multi,zhang2023corre}.

To illustrate the idea, let us consider the DNN approximation of a given band-limited target function $f(\bm x)$, $\bm x\in\mathbb R^d$, whose Fourier transform
\begin{equation}\label{fouriertransformdef}
	\hat f(\bm\xi):=\mathcal F[f](\bm\xi)=\int_{\mathbb R^d}f(\bm x)e^{-\ri 2\pi\bm\xi^{\rm T} \bm x}{\rm d}\bm x,
\end{equation}
has a compact support, i.e.,
\begin{equation}
	{\rm supp} \hat f(\bm{\xi})\subset B_{K}(\bm 0):=\{\bm{\xi}\in\mathbb R^d, |\bm{\xi}|\leq K\}.
\end{equation}
Note that the hyper-sphere $B_{K}(\bm 0)$ in the frequency domain can be partitioned into a union of $s+1$ concentric annulus with uniform or non-uniform radial dimension, e.g., for the case of uniform radial dimension,
\begin{equation}
	B_{K}(\bm 0)=\bigcup\limits_{j=0}^s A_j,\quad A_j:=\Big\{\bm{\xi}\in\mathbb R^d, \frac{jK}{s+1}\leq|\bm{\xi}|< \frac{(j+1)K}{s+1}\Big\}.
\end{equation}
Then, the target function in the frequency domain has a decomposition
\begin{equation}\label{decompositionfreq}
	\hat f(\bm{\xi})=\sum\limits_{j=0}^sI_{A_{j}}(\bm{\xi})\hat f(\bm{\xi}):=\sum\limits_{j=0}^s\hat f_j(\bm{\xi}),
\end{equation}
where $I_{A_j}(\bm{\xi})$ is the indicator function of the set $A_j$. From its definition,  the component $\hat f_j(\bm{\xi})$ has a ${\rm supp} \hat f_j(\bm{\xi})\subset A_j$, for $j=0, 1, \cdots, s$. A corresponding decomposition of \eqref{decompositionfreq} in the physical domain is given by
\begin{equation}\label{decompositionphysical}
	f(\bm x)=\sum\limits_{j=0}^sf_j(\bm x),
\end{equation}
with $f_j(\bm x)$ being the inverse Fourier transform
\begin{equation}
	f_j(\bm x)=\mathcal F^{-1}[\hat f_j](\bm x):=\int_{\mathbb R^d}\hat f_j(\bm \xi)e^{\ri 2\pi\bm\xi^{\rm T} \bm x}{\rm d}\bm \xi.
\end{equation}

With the decomposition \eqref{decompositionfreq}, an appropriate scaling can be used to transform the component $\hat f_j(\bm \xi)$ from  the high frequency  region $A_j$ to a low frequency region $A_j/\alpha_j$. The scaled version of $\hat f_j(\bm{\xi})$ is defined as
\begin{equation}\label{scalingfreq}
	\hat f_j^{\rm (scale)}(\bm \xi)=\hat f_j(\alpha_j\bm \xi),
\end{equation}
where $\alpha_j>1$ is an appropriate scaling factor for $A_j$. By the identity
\begin{equation}
	\mathcal F[g(a\bm x)](\bm\xi)=\Big(\frac{1}{|a|}\Big)^d\mathcal F[g]\Big(\frac{\bm\xi}{a}\Big),
\end{equation}
the scaling \eqref{scalingfreq} in the frequency domain leads to
\begin{equation}
	f^{\rm(scale)}_j(\bm x):=\mathcal F^{-1}[\hat f_j^{\rm (scale)}](\bm x)=\frac{1}{\alpha_j^d}f_j\Big(\frac{\bm x}{\alpha_j}\Big),
\end{equation}
or equivalently
\begin{equation}\label{physicalscale}
	f_j(\bm x)=\alpha_j^df^{\rm(scale)}_j(\alpha_j\bm x).
\end{equation}
By choosing an appropriate scale $\alpha_j$, we are able to make the Fourier spectrum of $\hat f_j^{\rm scale}(\bm \xi)$ into a lower frequency range, i.e.,
\begin{equation}
	{\rm supp}\hat f_j^{\rm scale}(\bm \xi)\subset \Big\{\bm{\xi}\in\mathbb R^d, \frac{jK}{(s+1)\alpha_j}\leq|\bm{\xi}|< \frac{(j+1)K}{(s+1)\alpha_j}\Big\}.
\end{equation}
\begin{figure}[htbp]
	\centering
	\includegraphics[width=\textwidth]{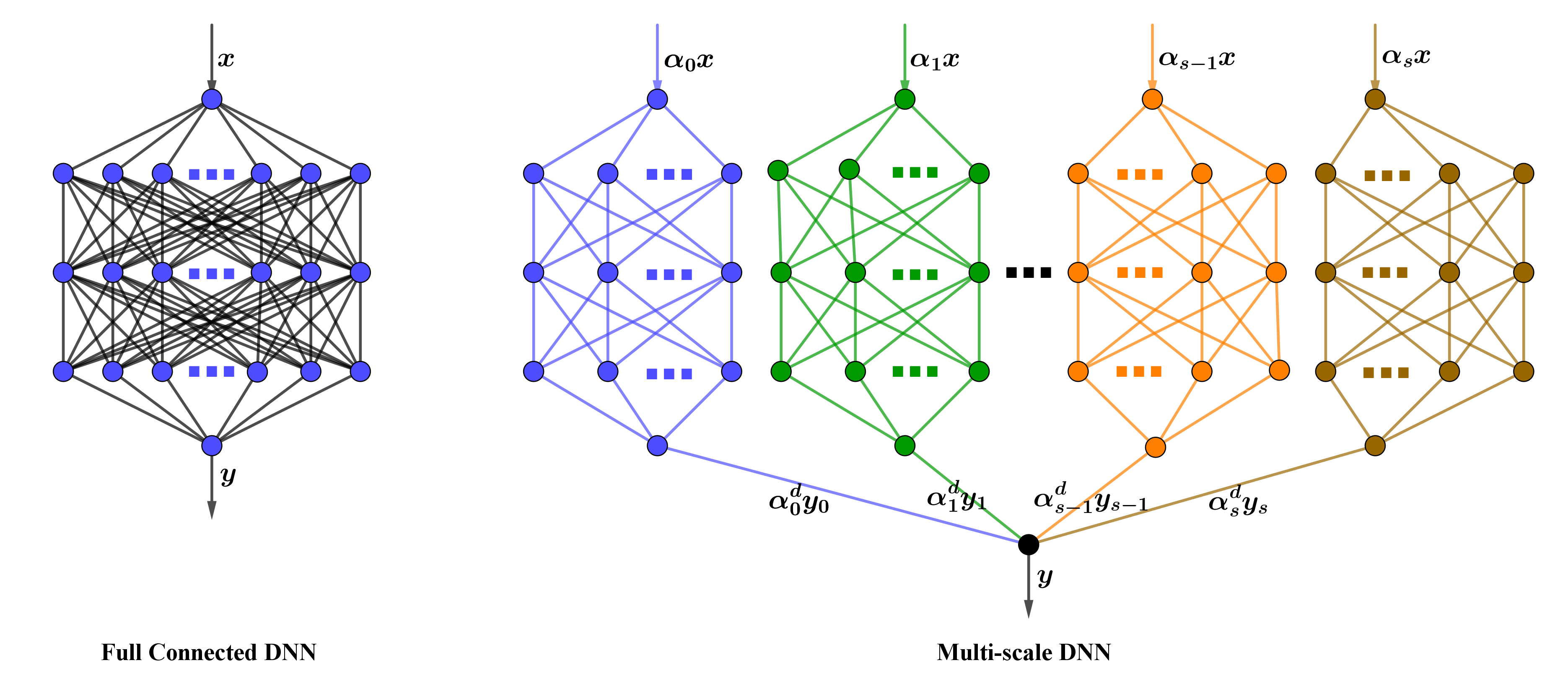}
	\caption{Sketches of a fully connected DNN and a Multi-scale DNN.}
	\label{net}
\end{figure}

As a result of the spectral bias of DNN, a fully connected DNN $f(\bm x;\bm\theta_{j})$ with parameters $\bm\theta_j$ can be trained to learn $f^{\rm(scale)}_j(\bm x)$ very fast if $(j+1)K/((s+1)\alpha_j)$ is small enough. Therefore, the decomposition \eqref{decompositionphysical} and scaling formula \eqref{physicalscale} implies that a deep learning algorithm using a neural network in the form
\begin{equation}\label{MscaleDNN}
	\mathcal N_s(\bm x;\bm{\theta})=\sum\limits_{j=0}^s\alpha_j^df(\alpha_j\bm x;\bm\theta_{j})
\end{equation}
can be expected to have a more uniform convergence and less spectral bias, i.e., frequency uniform  approximation to any band-limited function $f(\bm x)$. Deep neural networks defined by   \eqref{MscaleDNN} are named as the MscaleDNN and a schematic comparison between fully connected DNN and MscaleDNN is shown in Fig. \ref{net}.

Previous work presented in \cite{caixu2019,liu2020multi,wang2020multi} have shown that the MscaleDNN can reduce spectral bias significantly in learning highly oscillatory functions, however, mathematical analysis on the mechanism has not been presented in the literature. The following analysis will build a foundation for the MscaleDNN and provide a strategy to manipulate the neural networks.

\section{Error diffusion equation model of a two-layer MscaleDNN}
In this section,  the convergence of a machine learning algorithm for $d$-dimensional regression problems with two layers multi-scale neural networks is analyzed. We will show that the evolution of the error can be modeled by a diffusion equation in the Fourier frequency domain as the width of the network goes to infinity and learning rate approaches to zero.

Consider a regression problem with an objective function $y=f(\bm x)$ defined in a bounded domain $\Omega\subset\mathbb R^d$. The machine learning algorithm with a neural network denoted by $\mathcal N(\bm x,\bm{\theta})$ and mean square loss
\begin{equation}\label{meansquareloss}
	L(\bm\theta)=\frac{1}{2}\int_{\Omega}|\mathcal N(\bm x,{\bm \theta})-f(\bm x)|^2d\bm x,
\end{equation}
will be discussed in the following analysis.

The gradient descent dynamics based on the loss functional \eqref{meansquareloss} is
\begin{equation}\label{graddescentscheme}
	\bm\theta^{(k+1)}=\bm\theta^{(k)}-\tau\nabla L(\bm\theta^{(k)}),
\end{equation}
where $\tau$ is the learning rate. By regarding $\tau$ as the time step size, the continuum limit dynamics at $\tau\rightarrow 0$ is
\begin{equation}\label{graddescentdynamics}
	\frac{{\rm d}\bm\theta(t)}{{\rm d} t}=-\nabla L(\bm\theta(t)).
\end{equation}
With the mean square loss function \eqref{meansquareloss} and the chain rule of differentiation, we obtain
\begin{equation}\label{DNNfundynamics}
	\begin{split}
		\partial_t\mathcal N(\bm x,\bm\theta)=&[\nabla_{\bm\theta}\mathcal N(\bm x,\bm\theta)]^{\rm T} \frac{{\rm d}\bm\theta}{{\rm d} t}\\
		=&-\int_{\Omega}(\nabla_{\theta}\mathcal N(\bm x,\bm\theta))^{\rm T}\nabla_{\bm \theta}\mathcal N(\bm x',\theta)(\mathcal N(\bm x',\theta)-f(\bm x'))d\bm x'\\
		:=&-\int_{\Omega}\Theta(\bm x, \bm x';\bm\theta)(\mathcal N(\bm x',\bm \theta)-f(\bm x'))d\bm x',
	\end{split}
\end{equation}
for the dynamics of the network function $\mathcal N(\bm x,\theta)$, where
\begin{equation}\label{NTKdef}
	\Theta(\bm x, \bm x';\bm\theta)=(\nabla_{\theta}\mathcal N(\bm x,\bm\theta))^{\rm T}\nabla_{\bm \theta}\mathcal N(\bm x',\theta),
\end{equation}
is the neural tangent kernel (NTK) proposed in \cite{jacot2018neural}.

\begin{figure}[ht!]
	\center
	\subfloat{\includegraphics[scale=0.27]{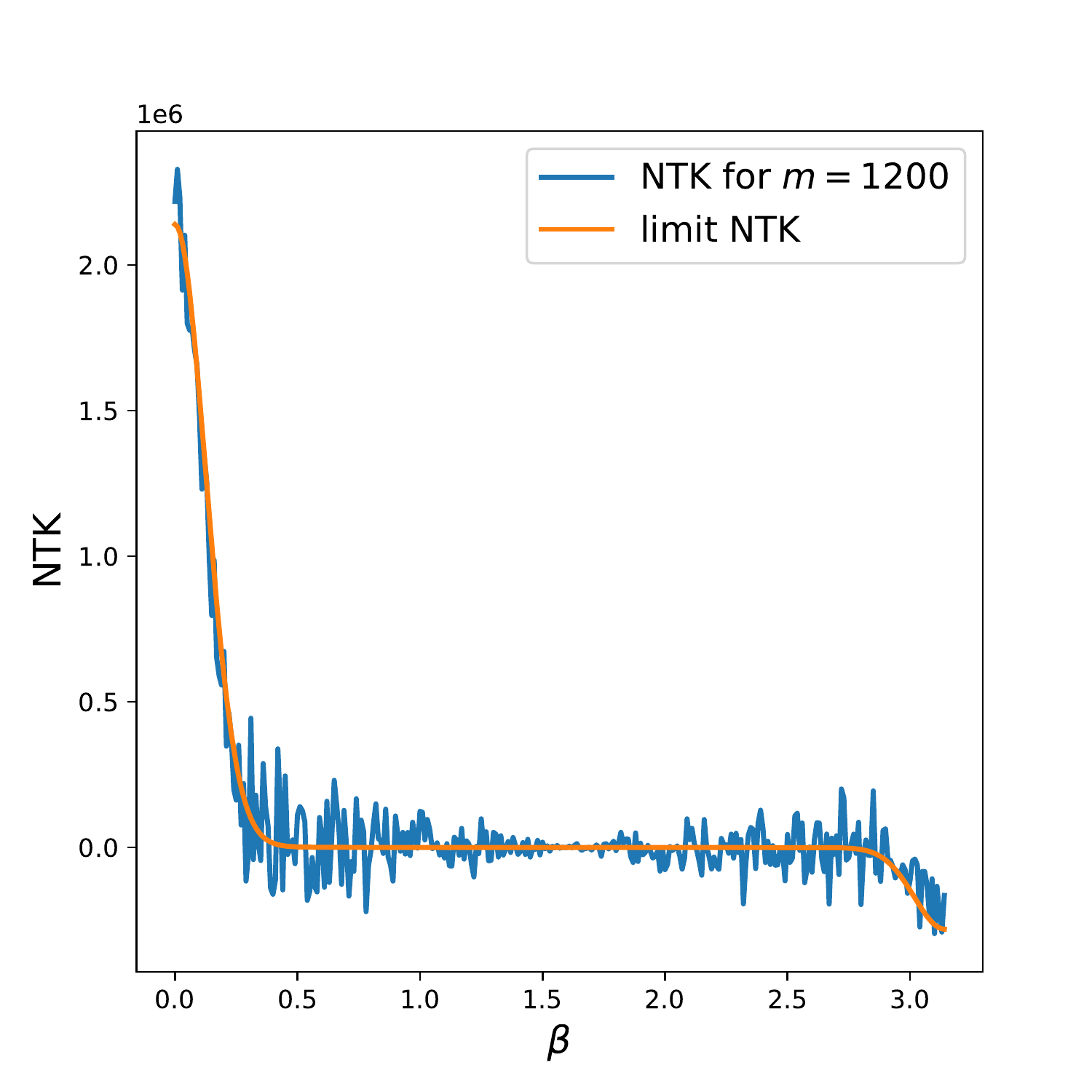}}
	\subfloat{\includegraphics[scale=0.27]{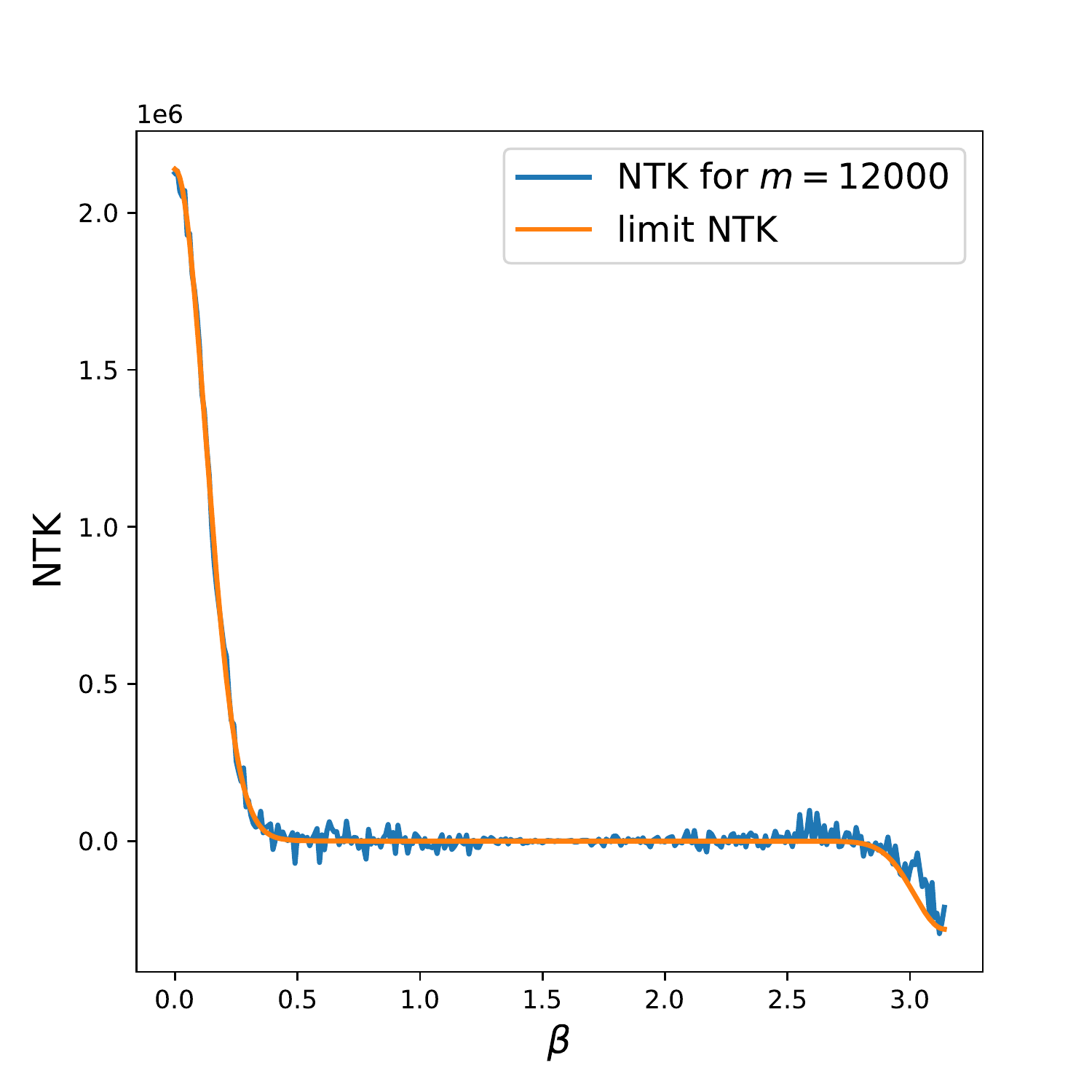}}
	\subfloat{\includegraphics[scale=0.27]{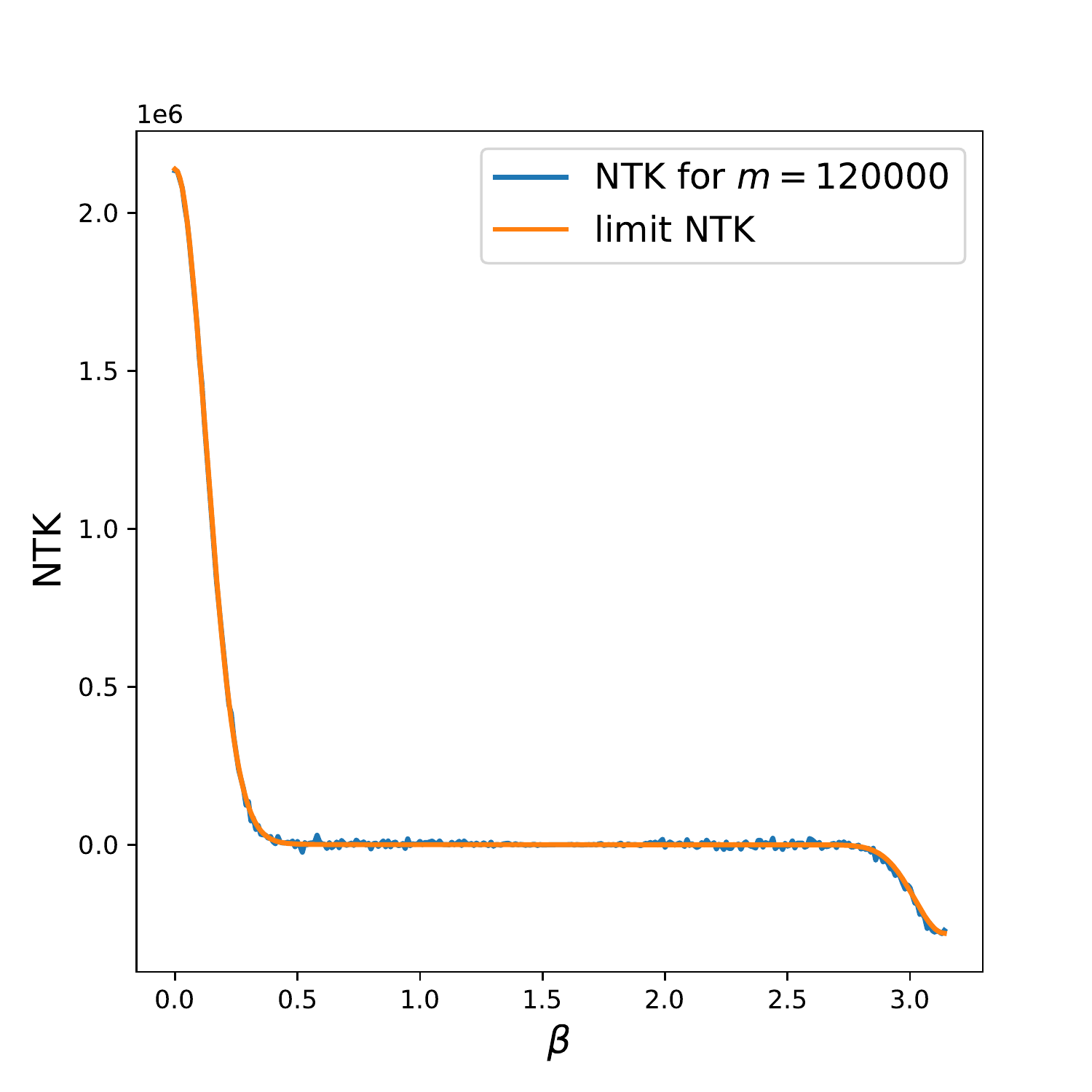}}
	\caption{The NTKs in \eqref{MsDNNNTK} and the limit NTK in \eqref{LIMITNTK} ($d=3, s=3$).}%
	\label{approxNTK}%
\end{figure}
\begin{figure}[ht!]
	\center
	\subfloat{\includegraphics[scale=0.27]{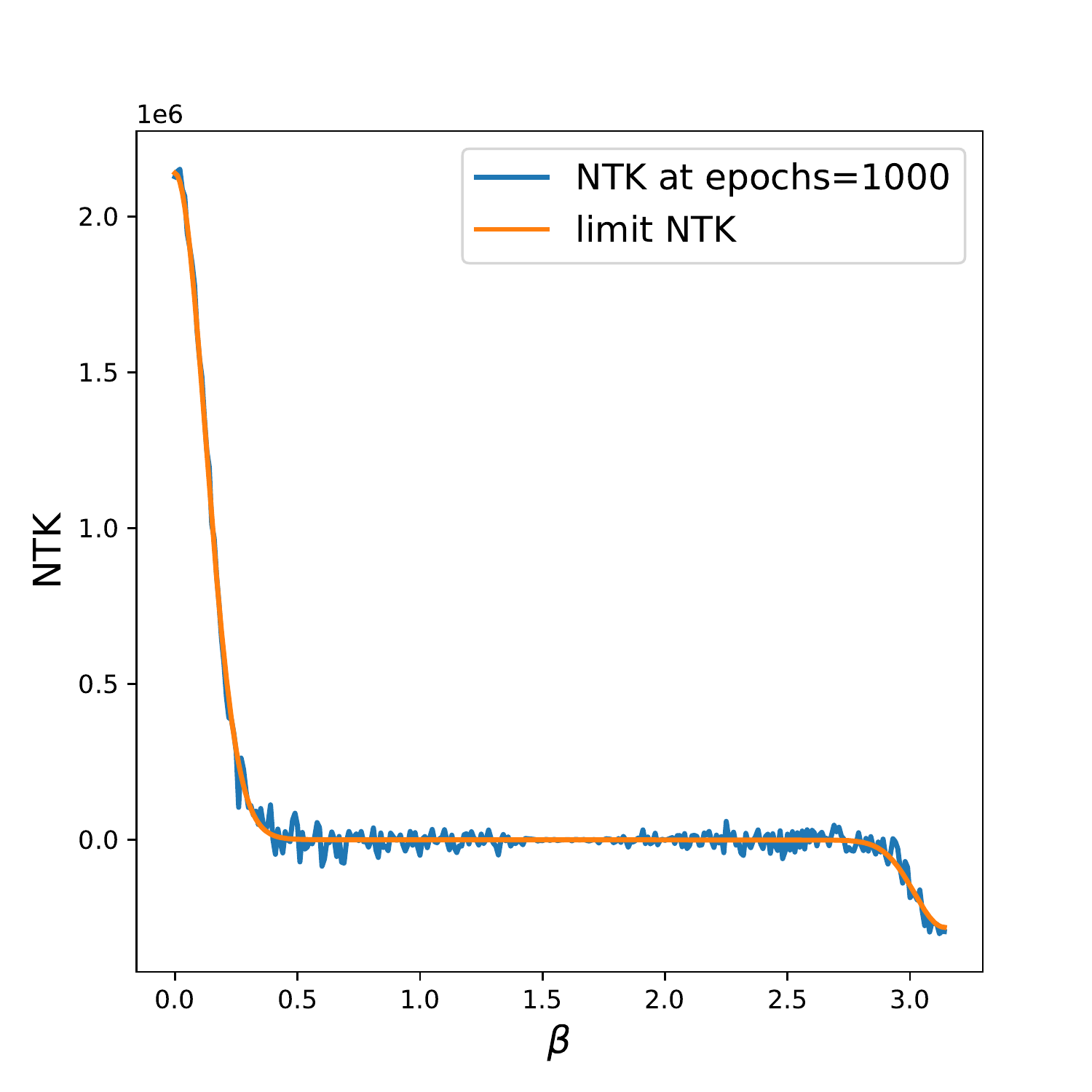}}
	\subfloat{\includegraphics[scale=0.27]{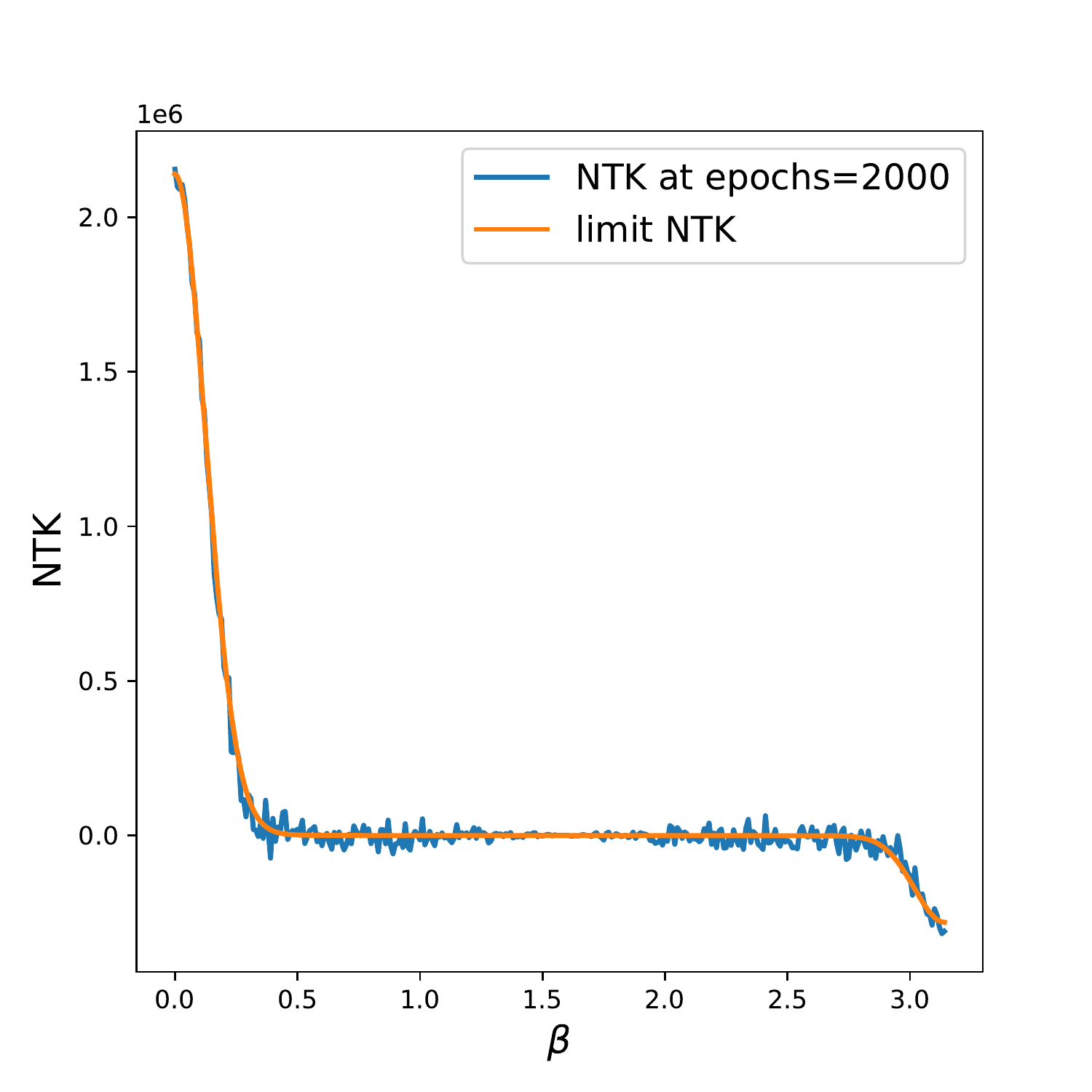}}
	\subfloat{\includegraphics[scale=0.27]{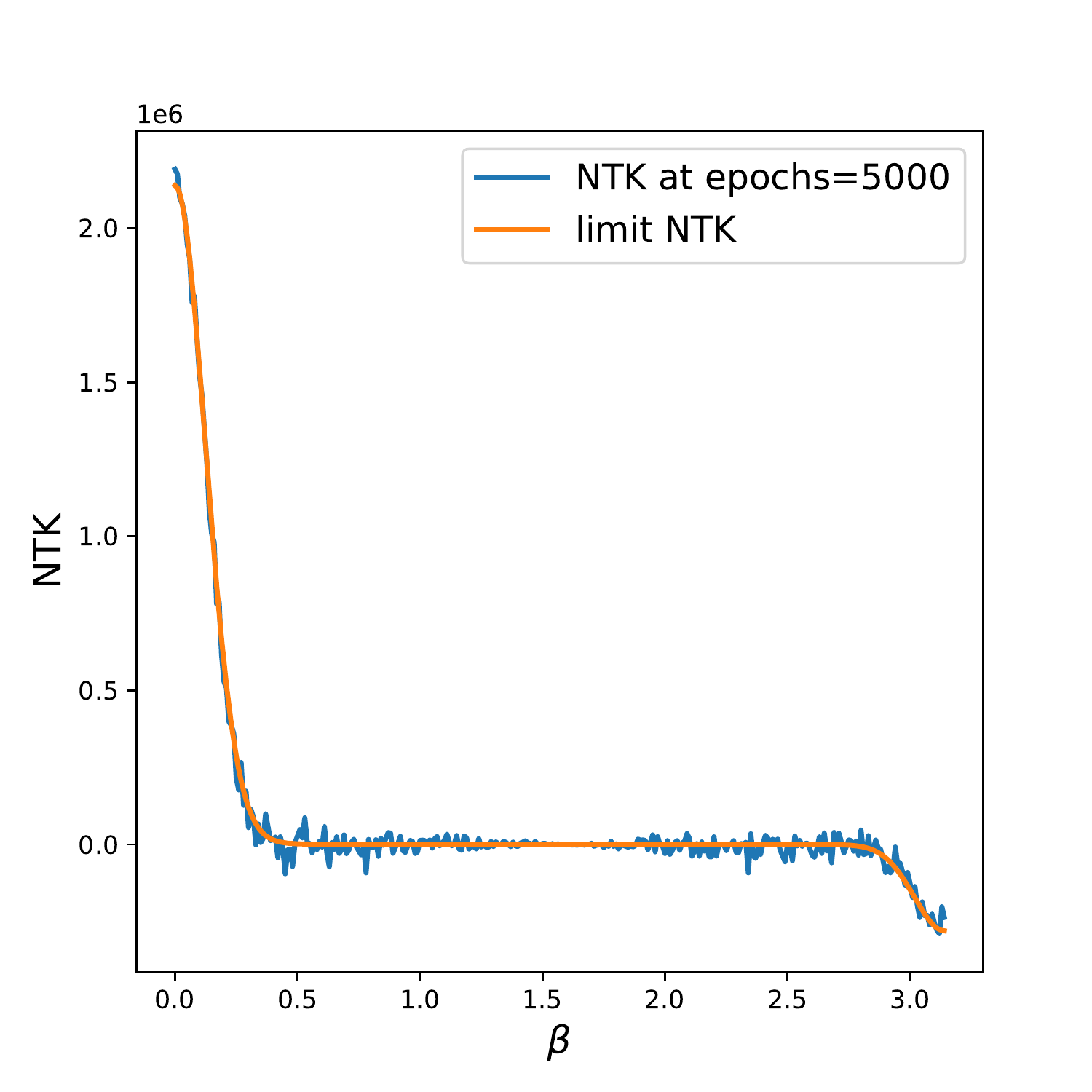}}
	\caption{The NTKs in \eqref{MsDNNNTK}  with $d=3, s=3, N=12000$ during training.}%
	\label{NTKstatic}%
\end{figure}

A multi-scale neural network with one hidden layer (see. Fig. \ref{net} (right)) 
is given as
\begin{equation}\label{MsDNNfun}
	{\mathcal N}_s(\bm x,\bm\theta)=\frac{1}{\sqrt{N}}\sum\limits_{j=0}^{s}\alpha_j^{d}\sum\limits_{k=1}^{q}\sigma(\bm\theta_{jq+k}^{\rm T}\alpha_j\bm x+b_{jq+k}),\quad \bm x\in \Omega:=[-1, 1]^d,
\end{equation}
where $s+1$ is the number of scales, $\{\alpha_j\}_{j=0}^s$ are the scaling factors, $q$ is the number of neurons for each scale, $N=(s+1)q$ is the total number of neurons in the hidden layer. Apparently, the network includes the standard fully connected neural network with one hidden layer as a special case of $s=0$.
For this two-layer multi-scale neural network, a direct calculation gives its NTK
\begin{equation}\label{MsDNNNTK}
	{\Theta}_s(\bm x, \bm x';\bm\theta)=\sum\limits_{j=0}^{s}\frac{\alpha_j^{2d}(\alpha_j^{2}\bm x^{\rm T}\bm x'+1)}{N}\sum\limits_{k=1}^{q}\sigma'(\bm\theta_{jq+k}^{\rm T}\alpha_p\bm x+b_{jq+k})\sigma'(\bm\theta_{jq+k}^{\rm T}\alpha_j\bm x'+b_{jq+k}).
\end{equation}
Setting the activation function
\begin{equation}
   \sigma(x)=\sin(x) \label{sine_act}
\end{equation}
and assuming all the parameters $\{\theta_{jk}\}$ in $\bm\theta_j=(\theta_{j1}, \theta_{j2},\cdots, \theta_{jd})^{\rm T}$, $\{b_j\}$ are independent random variables of normal distribution,  then, by the law of large numbers and identity
$$(2\pi)^{-\frac{d+1}{2}}\int_{\mathbb R^{d+1}}e^{\ri(\bs\theta^{\rm T}\bs x+yb)}e^{-\frac{|\bs \theta|^2+b^2}{2}}d\bs\theta db=e^{-\frac{|\bs x|^2+y^2}{2}},\quad \forall \bs x\in\mathbb R^d,\;\; y\in \mathbb R,$$
we have
\begin{equation}\label{LIMITNTK}
	\begin{split}
		\lim\limits_{q\rightarrow\infty}{\Theta}_s(\bm x, \bm x';\bm\theta)
		=&\sum\limits_{j=0}^{s}\frac{\alpha_j^{2d}(\alpha^2_{j}\bm x^{\rm T}\bm x'+1)}{s+1}\mathbb E(\cos(\bm\theta_{1}^{\rm T}\alpha_j\bm x+b_{1})\cos(\bm\theta_{1}^{\rm T}\alpha_j\bm x'+b_{1}))\\
		=&\sum\limits_{j=0}^{s}\frac{\alpha_j^{2d}(\alpha^2_{j}\bm x^{\rm T}\bm x'+1)}{2(s+1)}\Big[e^{-2}e^{-\frac{\alpha^2_{j}|\bm x+\bm x'|^2}{2}}+e^{-\frac{\alpha^2_{j}|\bm x-\bm x'|^2}{2}}\Big].
	\end{split}
\end{equation}
Apparently, $\{\bs\theta_1, b_1\}$ can replaced by the parameters of any neuron in the hidden layer.
According to the analysis in \cite{jacot2018neural}, the NTK will be static during the training assuming the width of the neural network tends to infinity. In addition, the limit NTK is also a convolution kernel as presented in \cite{cho2009kernel, ronen2019convergence}. Suppose $\bs x$ and $\bs x'$ are located on the unit sphere, i.e., $|\bs x|=|\bs x'|=1$, then the limit NTK is a function of the angle $\beta$ between $\bs x$ and $\bs x'$. The NTKs of some multi-scale neural networks with finite width are compared with their infinite width limit in Fig. \ref{approxNTK}. We can see that the NTK \eqref{MsDNNNTK} has an  limit given above as $q\rightarrow\infty$.
In order to validate the static property of the limit NTK, we train a multi-scale neural network with $s=3$, $N=12000$ to fit a 3-dimensional function in the domain $[-1, 1]^3$. The scaling parameters $\alpha_p$ are set to be $2^p$. The NTKs of the multi-scale neural network after training $1000, 2000, 5000$ epochs are compared with the limit NTK in Fig. \ref{NTKstatic}. The results clearly show that the NTK is static during training.

Consequently, as the width of the network goes to infinity, the dynamics of the gradient descent learning \eqref{DNNfundynamics} tends to
\begin{equation}\label{limitdynamicMsDNNbiastrue}
	\begin{split}
		&\partial_t(\mathcal N_s(\bm x,\bm\theta)-f(\bm x))\\
		=&-\sum\limits_{j=0}^{s}\frac{\alpha_j^{2(d+1)}\bm x^{\rm T}}{2(s+1)}\int_{\Omega}\Big[e^{-2}\mathcal G_j(\bm x+\bm x')+\mathcal G_j(\bm x-\bm x')\Big]\bm x'(\mathcal N_s(\bm x',\bm\theta)-f(\bm x')){\rm d}\bm x'\\
		&-\sum\limits_{j=0}^{s}\frac{\alpha_j^{2d}}{2(s+1)}\int_{\Omega}\Big[e^{-2}\mathcal G_j(\bm x+\bm x')+\mathcal G_j(\bm x-\bm x')\Big](\mathcal N_s(\bm x',\bm\theta)-f(\bm x')){\rm d}\bm x',
	\end{split}
\end{equation}
where
\begin{equation}\label{Gaussion}
	\mathcal G_j(\bm x):=e^{-\alpha^2_j|\bm x|^2/2},\quad \bm x\in\mathbb R^d,
\end{equation}
is the scaled Gaussian function.

Next, we define a zero extension of the error function by
\begin{equation}
	\eta(\bm x,\bm \theta)=\begin{cases}
		0, & \bm x\notin\Omega,\\
		\mathcal N_s(\bm x,\bm\theta)-f(\bm x), & \bm x\in\Omega,
	\end{cases}
\end{equation}
then, the dynamic system \eqref{limitdynamicMsDNNbiastrue} can be rewritten as
\begin{equation}\label{limitdynamicMsDNNbiastruere}
	\begin{split}
		& \partial_t\eta(\bm x,\bm\theta)I_{\Omega}(\bm x) \\
		=&-\sum\limits_{j=0}^{s}\frac{I_{\Omega}(\bm x)\alpha_j^{2(d+1)}\bm x^{\rm T}}{2(s+1)}\int_{\mathbb R^d}\Big[e^{-2}\mathcal G_j(\bm x+\bm x')+\mathcal G_j(\bm x-\bm x')\Big]\bm x'\eta(\bm x',\bm\theta){\rm d}\bm x'\\
		&-\sum\limits_{j=0}^{s}\frac{I_{\Omega}(\bm x)\alpha_j^{2d}}{2(s+1)}\int_{\mathbb R^d}\Big[e^{-2}\mathcal G_j(\bm x+\bm x')+\mathcal G_j(\bm x-\bm x')\Big]\eta(\bm x',\bm\theta){\rm d}\bm x',
	\end{split}
\end{equation}
where an indicator function $I_{\Omega}(\bm x)$ is used to extend the equation to the whole space.

Existing works on DNN convergence analysis employ a discrete version  of \eqref{limitdynamicMsDNNbiastruere} in the physical space by analyzing the eigenvalues of the Gram matrix \cite{tancik2020fourier,luo2022exact,peng2023non,ronen2019convergence, lee2019wide}.
However, to get a precise information on the spectral bias phenomena, it is more natural to study the convergence behavior in the Fourier domain as follows.

Given any $g(\bm x)\in L^1(\mathbb R^d)$, the Fourier transform defined in \eqref{fouriertransformdef} has the following identities
\begin{equation}\label{derivativefourier}
	\mathcal F[\nabla g](\bm\xi)=2\pi\ri\bm\xi\mathcal F[g](\bm\xi),\quad \nabla \hat g(\bm\xi)=-2\pi\ri \mathcal F[\bm xg(\bm x)](\bm\xi),\quad \forall \bm xg(\bm x)\in (L^1(\mathbb R^d))^d,
\end{equation}
and
\begin{equation}\label{gaussianfourier}
	\mathcal F[e^{-|\bm x|^2}](\xi)=\pi^{\frac{d}{2}}e^{-\pi^2|\bm\xi|^2},\quad \mathcal F[g(a\bm x)](\bm\xi)=\Big(\frac{1}{|a|}\Big)^d\mathcal F[g]\Big(\frac{\bm\xi}{a}\Big).
\end{equation}
In addition, given two functions $h(\bm x)$, $g(\bm x)$, their cross-correlation and convolution are defined as
\begin{equation}
	h\star g:=\int_{\mathbb R^d}\overline{h(\bm x')}g(\bm x+\bm x'){\rm d}\bm x', \quad h* g:=\int_{\mathbb R^d}h(\bm x')g(\bm x-\bm x'){\rm d}\bm x',
\end{equation}
and we have teh following identities,
\begin{equation}\label{ccorrconvtheorem}
	\widehat{h\star g}(\bm \xi)=\overline{\hat h(\bm\xi)}\hat g(\bm\xi), \quad  \widehat{h* g}(\bm\xi)=\hat h(\bm\xi)\hat g(\bm\xi),\quad \widehat{fg}(\bm{\xi})=f*g(\bm \xi).
\end{equation}
Taking Fourier transform \eqref{fouriertransformdef} on both sides of \eqref{limitdynamicMsDNNbiastruere} with respect to $\bm x$ and then applying \eqref{derivativefourier}-\eqref{ccorrconvtheorem} to rearrange the terms gives a integral-differential equation
\begin{equation}\label{MsDNNfreqdynamicsbaistrue}
	\begin{split}
		\frac{\partial \hat{\eta}(\bm\xi,\bm\theta(t))*\hat I_{\Omega}(\bm{\xi})}{\partial t}=&\left[\nabla_{\bm\xi}\cdot\Big[\sum\limits_{j=0}^{s}\frac{\alpha_j^{2(d+1)}\widehat{\mathcal G}_j(\bm \xi)}{8\pi^2(s+1)}\Big(\nabla_{\bm \xi}\hat{ \eta}(\bm\xi,\bm\theta(t))-e^{-2}\nabla_{\bm \xi}\overline{\hat{ \eta}(\bm\xi,\bm\theta(t))}\Big)\Big]\right.\\
		&\left.-\sum\limits_{j=0}^{s}\frac{\alpha_j^{2d}\widehat{\mathcal G}_j(\bm \xi)}{2(s+1)}[e^{-2}\bar{\hat{ \eta}}(\bm\xi,\bm\theta(t))+\hat{ \eta}(\bm\xi,\bm\theta(t))]\right]*\hat I_{\Omega}(\bm{\xi}),
	\end{split}
\end{equation}
where
\begin{equation}\label{FourierGauss}
	\widehat{\mathcal G}_j(\bm\xi)=(2\pi)^{\frac{d}{2}}\alpha_j^{-d}e^{-\frac{2\pi^2|\bm \xi|^2}{\alpha^2_j}}.
\end{equation}

The convolution with $\hat I_{\Omega}(\bm \xi)$ makes the model too complicate to analyze. Nevertheless, if we consider the limit of infinite large domain, i.e., $\Omega\rightarrow \mathbb R^d$, the limit of $\hat I_{\Omega}(\bm \xi)$ is the Dirac delta function $\delta(\xi)$ and then \eqref{MsDNNfreqdynamicsbaistrue} simplifies to
\begin{equation}\label{MsDNNfourierfreqdynamicsbaistrue}
	\begin{split}
		\frac{\partial \hat{\eta}(\bm\xi,\bm\theta(t))}{\partial t}=&\nabla_{\bm\xi}\cdot\Big[\sum\limits_{j=0}^{s}\frac{\alpha_j^{2(d+1)}\widehat{\mathcal G}_j(\bm \xi)}{8\pi^2(s+1)}\Big(\nabla_{\bm \xi}\hat{ \eta}(\bm\xi,\bm\theta(t))-e^{-2}\nabla_{\bm \xi}\overline{\hat{ \eta}(\bm\xi,\bm\theta(t))}\Big)\Big]\\
		&-\sum\limits_{j=0}^{s}\frac{\alpha_j^{2d}\widehat{\mathcal G}_j(\bm \xi)}{2(s+1)}[e^{-2}\bar{\hat{ \eta}}(\bm\xi,\bm\theta(t))+\hat{ \eta}(\bm\xi,\bm\theta(t))].
	\end{split}
\end{equation}
Define
\begin{equation}\label{DNNcoefficients}
	\begin{split}
		A_s^{\pm}(\bm\xi)=\frac{1\pm e^{-2}}{8\pi^2(s+1)}\sum\limits_{j=0}^{s}\alpha_j^{2(d+1)}\widehat{\mathcal G}_j(\bm\xi),\quad  B_s^{\pm}(\bm\xi)=\frac{1\pm e^{-2}}{2(s+1)}\sum\limits_{j=0}^{s}\alpha_j^{2d}\widehat{\mathcal G}_j(\bm\xi),
	\end{split}
\end{equation}
and denote by $\hat\eta^{\pm}(\bm\xi,\bm\theta(t))$ the real and imaginary parts of $\hat{\eta}(\bm\xi,\bm\theta(t))$, i.e.,
$$\hat{\eta}(\bm\xi,\bm\theta(t))=\hat{\eta}^+(\bm\xi,\bm\theta(t))+\ri \hat{\eta}^-(\bm\xi,\bm\theta(t)).$$

\medskip

\noindent ({\bf  Diffusion Model}) As the coefficients in \eqref{MsDNNfourierfreqdynamicsbaistrue} are real valued functions, we can rewrite \eqref{MsDNNfourierfreqdynamicsbaistrue} into two independent equations
\begin{equation}\label{realimagdynamicsbaistrueMsDNN}
	\partial_t \hat\eta^{\pm}(\bm\xi,t)=\nabla_{\bm\xi}\cdot\Big[A_s^{\mp}(\bm\xi)\nabla_{\bm\xi} \hat\eta^{\pm}(\bm\xi,t)\Big]-B_s^{\pm}(\bm\xi)\hat\eta^{\pm}(\bm\xi,t),\quad \bm\xi\in\mathbb R^d,
\end{equation}
with respect to the real and imaginary parts of $\hat{\eta}(\bm\xi,\bm\theta(t))$, respectively.

A simpler diffusion equation can be  derived if the bias are set to zero in the network. In fact, a function represented by the network without bias has the form
\begin{equation}\label{msdnnwithoutbias}
	{\mathcal N}_s(\bm x,\bm\theta)=\frac{1}{\sqrt{N}}\sum\limits_{j=0}^{s}\alpha_j^{d}\sum\limits_{k=1}^{q}\sigma(\bm\theta_{jq+k}^{\rm T}\alpha_j\bm x),\quad \bm x\in \Omega:=[-1, 1]^d,
\end{equation}
and the neural tangent kernel is given by
\begin{equation}
	{\Theta}_s(\bm x, \bm x';\bm\theta)=\frac{\bm x^{\rm T} \bm x'}{N}\sum\limits_{j=0}^{s}\alpha_j^{2(d+1)}\sum\limits_{k=1}^{q}\sigma'(\bm\theta_{jq+k}^{\rm T}\alpha_j\bm x)\sigma'(\bm\theta_{jq+k}^{\rm T}\alpha_j\bm x').
\end{equation}
Setting the activation function $\sigma(x)=\sin(x)$ again, and assuming all the parameters $\{\theta_p\}$ are independent random variables of normal distribution, then, by law of large numbers, we have
\begin{equation}
	\begin{split}
		\lim\limits_{q\rightarrow\infty}{\Theta}_s(\bm x, \bm x';\bm\theta)&=\lim\limits_{q\rightarrow\infty}\frac{\bm x^{\rm T} \bm x'}{N}\sum\limits_{j=0}^{s}\alpha_j^{2(d+1)}\sum\limits_{k=1}^{q}\cos(\bm\theta_{jq+k}^{\rm T}\alpha_j\bm x)\cos(\bm\theta_{jq+k}^{\rm T}\alpha_j\bm x')\\
		&=\frac{\bm x^{\rm T} \bm x'}{2(s+1)}\sum\limits_{j=0}^{s}\alpha_j^{2(d+1)}\mathbb E(\cos(\bm\theta_{1}^{\rm T}\alpha_j\bm x)\cos(\bm\theta_{1}^{\rm T}\alpha_j\bm x'))\\
		&=\frac{\bm x^{\rm T} \bm x'}{2(s+1)}\sum\limits_{j=0}^{s}\alpha_j^{2(d+1)}\Big[\mathcal G_j(\bm x+\bm x')+\mathcal G_j(\bm x-\bm x')\Big].
	\end{split}
\end{equation}
As the width of the network goes to infinity, the dynamics of the gradient descent learning tends to
\begin{equation}\label{limitdynamicMsDNN}
	\partial_t\eta(\bm x,\theta)
	=-\frac{\bm x^{\rm T}}{2(s+1)}\int_{\Omega}\sum\limits_{j=0}^{s}\alpha_j^{2(d+1))}\Big[\mathcal G_j(\bm x+\bm x')+\mathcal G_j(\bm x-\bm x')\Big]\bm x'\eta(\bm x',\bm\theta){\rm d}\bm x'.
\end{equation}
Mimicking the derivation for \eqref{MsDNNfourierfreqdynamicsbaistrue}, we obtain from \eqref{limitdynamicMsDNN} that
\begin{equation}
	\begin{split}
		\frac{\partial \hat{\eta}(\bm\xi,\bm\theta(t))}{\partial t}=&\nabla_{\bm\xi}\cdot\Big[\sum\limits_{j=0}^{s}\frac{\alpha_j^{2(d+1))}\widehat{\mathcal G}_j(\bm\xi)}{8\pi^2(s+1)}\Big(\nabla_{\bm\xi}\hat{ \eta}(\bm\xi,\bm\theta(t))-\nabla_{\bm\xi}\overline{\hat{\eta}(\bm\xi,\bm\theta(t))}\Big)\Big]\\
		=&\ri\nabla_{\bm\xi}\cdot\Big[\sum\limits_{j=0}^{s}\frac{ \alpha_j^{2(d+1))}}{4\pi^2(s+1)}\widehat{\mathcal G}_j(\bm\xi)\nabla_{\bm \xi}\hat{\eta}^-(\bm\xi,\bm\theta(t))\Big],
	\end{split}
\end{equation}
where $\hat{\eta}^-(\bm\xi,\bm\theta(t)):=\mathfrak{Im}\Big\{\hat{\eta}(\bm\xi,\bm\theta(t))\Big\}$. The dynamic system \eqref{limitdynamicMsDNN} in the Fourier frequency domain implies that only the imaginary part of the error evolves during the gradient descent training if a two layer multi-scale neural network with activation function $\sigma(x)=\sin(x)$ and zero bias is used. This conclusion is consistent with the fact that the network function \eqref{msdnnwithoutbias} can only be used to fit odd functions. Actually, the necessity of non-zero biases in a two layer neural network has been emphasized in \cite{ronen2019convergence, lee2019wide}.

\begin{figure}[ht!]
	\center
	\subfloat{\includegraphics[scale=0.21]{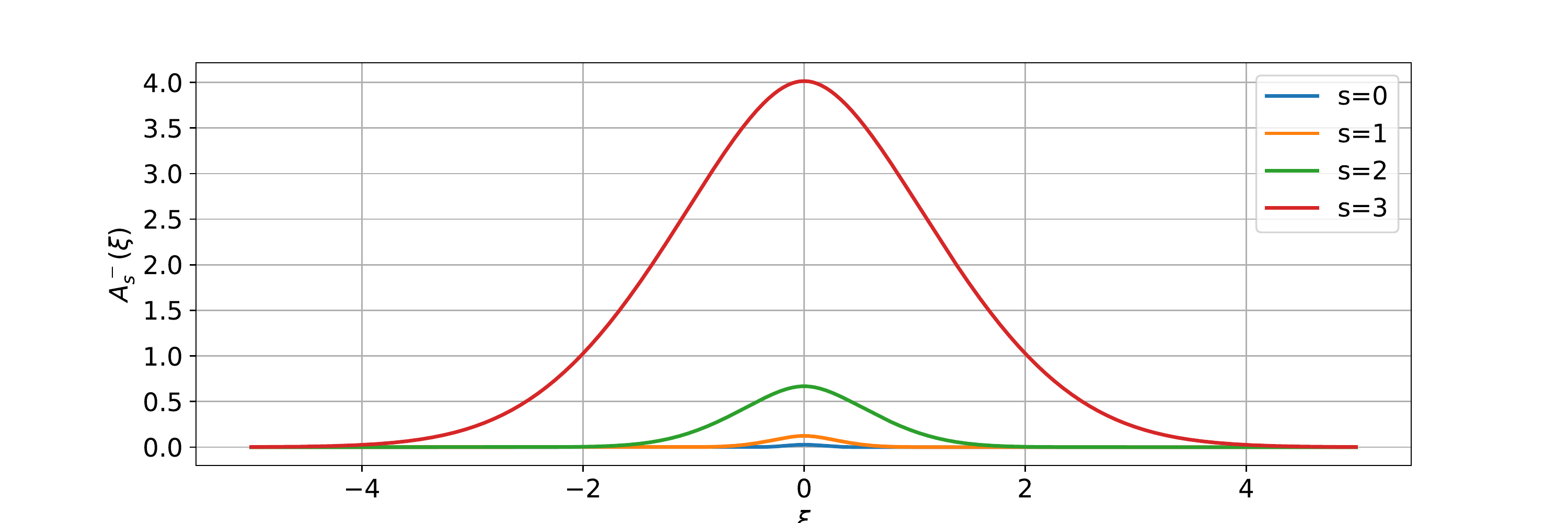}}
	\subfloat{\includegraphics[scale=0.21]{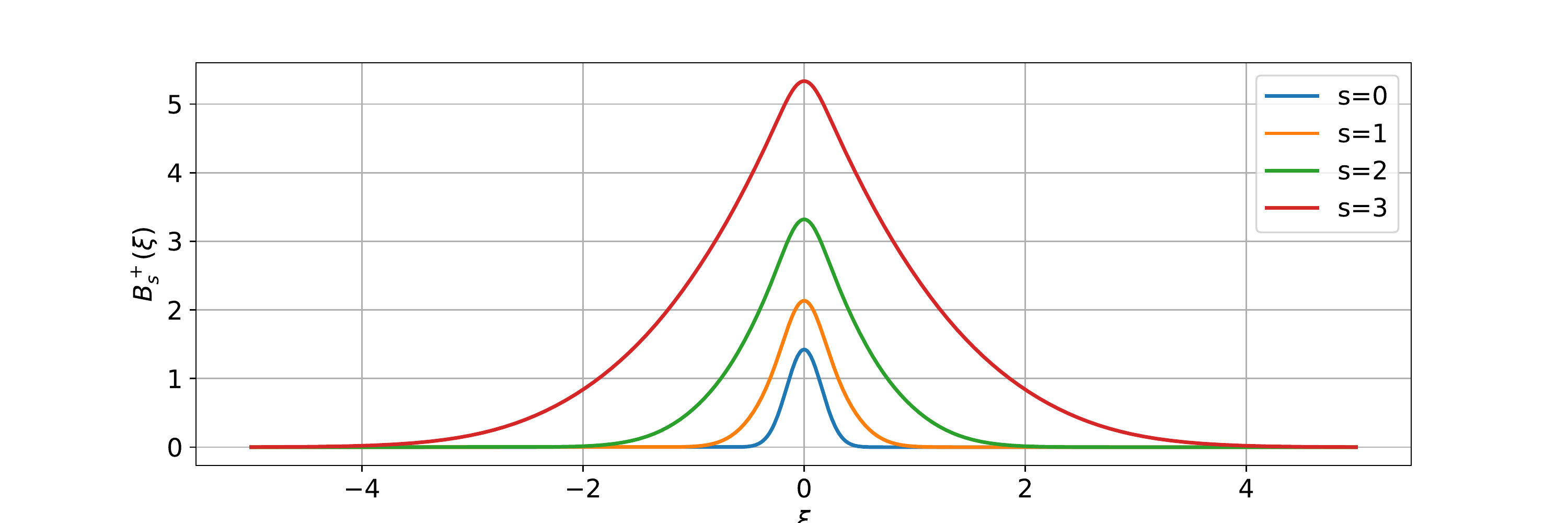}}\\
	\subfloat{\includegraphics[scale=0.21]{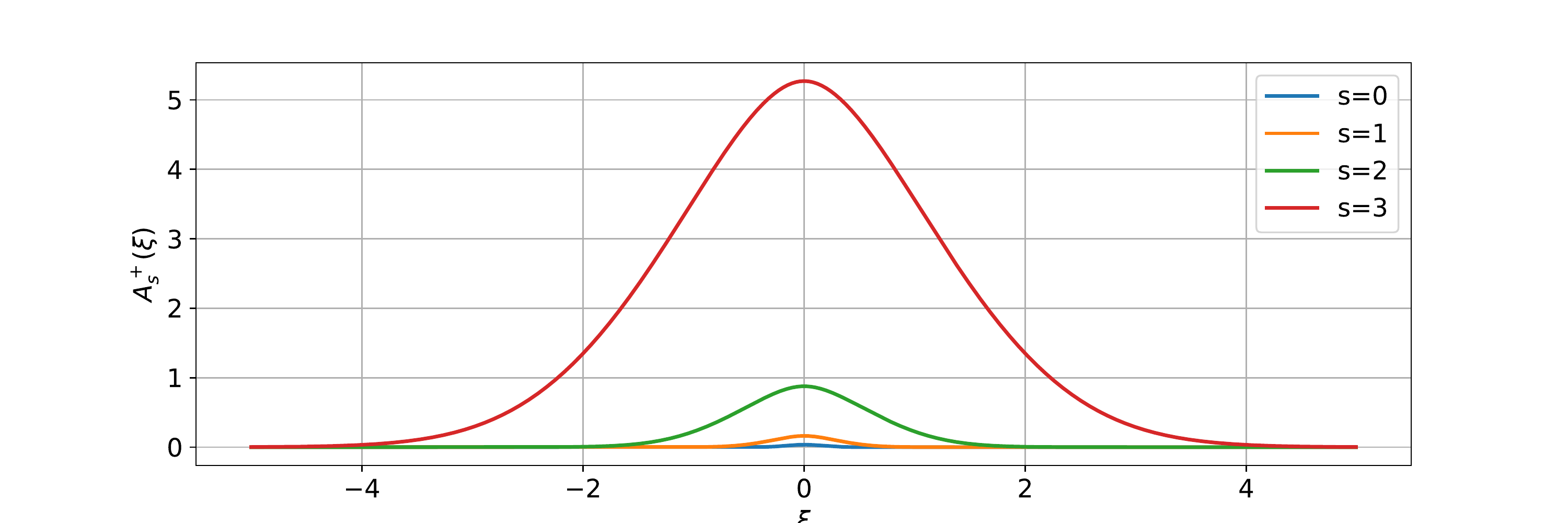}}
	\subfloat{\includegraphics[scale=0.21]{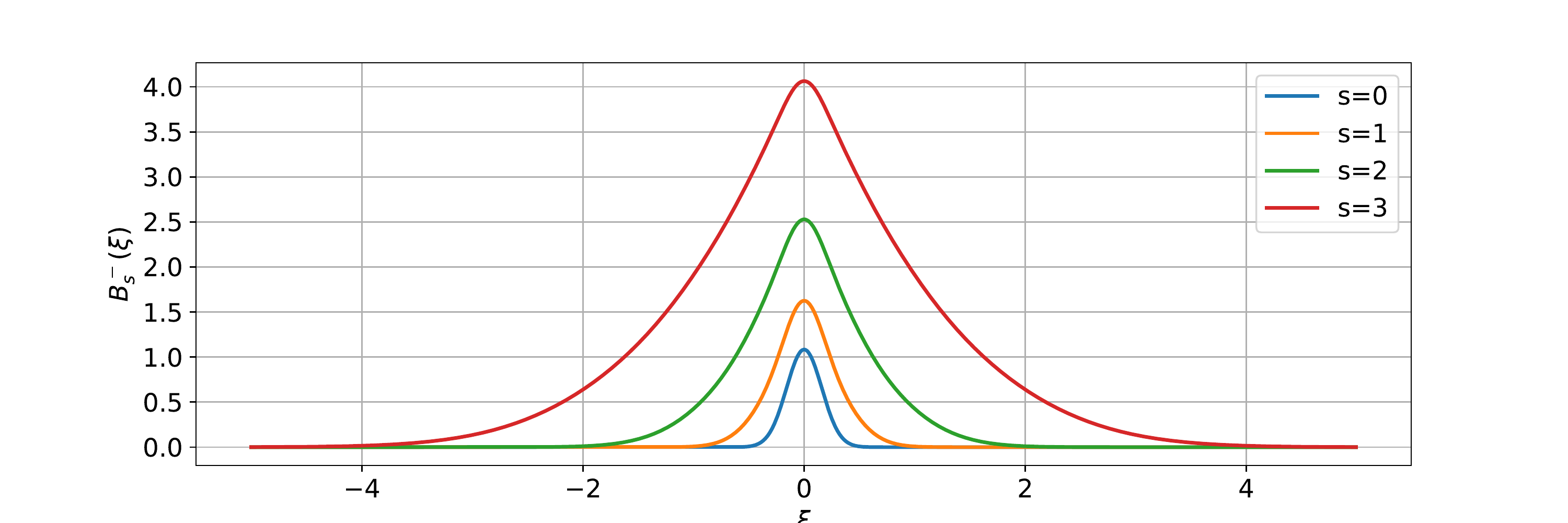}}
	\caption{Diffusion coefficients $A^{\mp}_s(\xi)$ (left) and $B^{\pm}_s(\xi)$ (right) with $\alpha_j=2^j$,  $s=0, 1, 2, 3$.}%
	\label{coefficientsplot}%
\end{figure}

Note that $A_s^{\pm}(\bm\xi), B^{\mp}_s(\bm\xi)$ defined in \eqref{DNNcoefficients} are positive functions in $\mathbb R^d$. Therefore, the solution of \eqref{realimagdynamicsbaistrueMsDNN} has an energy equality
\begin{equation}\label{energyequality}
	\frac{\rm d}{{\rm d} t}\int_{\mathbb R^d}|\hat\eta^{\pm}(\bm\xi, t)|^2{\rm d}\bm\xi=-2\int_{\mathbb R^d}\Big[A^{\mp}_s(\bm\xi)\Big|\nabla_{\bm\xi} \hat\eta^{\pm}(\bm\xi, t)\Big|^2+B^{\pm}_s(\bm\xi)|\hat\eta^{\pm}(\bm\xi, t)|^2\Big]{\rm d}\bm\xi,
\end{equation}
which implies that the solution $\hat\eta^{\pm}(\bm\xi, t)\rightarrow 0$ for any $\bm\xi\in\mathbb R^d$ as $t\rightarrow\infty$. That means the gradient descent learning for a fitting problem with one hidden layer neural network is convergent assuming that the learning rate is sufficiently small and the width of the neural network is sufficiently large. It is clear that the diffusion coefficients $\{A_s^{\pm}(\bm\xi), B_s^{\mp}(\bm\xi)\}$ plays a key role in the error decay speed. Several plots of the coefficients $\{A_s^{\mp}(\xi), B^{\pm}_s(\xi)\}$ are given in Fig. \ref{coefficientsplot} for different scales. We can see that both $A_s^{\mp}(\xi)$ and  $B^{\pm}_s(\xi)$ have larger support and maximum values with an increasing scale $s$. This implies that larger $s$ will leads to fast error reduction in a wider frequency region.

\section{Spectral bias analysis of a two layer MscaleDNN using the diffusion equation model}

The analysis in previous section has shown that the error dynamics of the gradient descent learning can be approximately described by the diffusion equations \eqref{realimagdynamicsbaistrueMsDNN} in the Fourier spectral domain when the network width and the domain size go to infinity and the learning rate to zero. In this section, we will first propose a Hermite spectral method to obtain highly accurate numerical solutions of the diffusion equation. Some numerical results will be presented to show that the error dynamics predicted by the diffusion model matches well with that of the MscaleDNN during realistic training. Moreover, the results also validate the capability of spectral bias reduction of the MscaleDNNs for wider range of frequencies.  For simplicity, we only consider the 1-dimensional case to illustrate the main results.


\subsection{Hermite spectral method for the  diffusion equation problem}
In order to examine quantitatively the decay of the error in the Fourier domain, we will solve numerically the equations in \eqref{realimagdynamicsbaistrueMsDNN} with a Hermite spectral method for the $\xi$-variable of the equations in \eqref{realimagdynamicsbaistrueMsDNN} on the unbounded computational domain. For this purpose, we introduce the Hermite functions (cf. \cite{ShenTaoWang2011}) defined by
\begin{equation}
	\widehat H_n(\xi)=\frac{1}{\pi^{1/4}\sqrt{2^nn!}}e^{-\xi^2/2}H_n(\xi),\quad n\geq 0, \;\; \xi\in\mathbb R,
\end{equation}
where $H_n(\xi)$ are Hermite polynomials. The Hermite functions $\widehat H_n(\xi)$ are orthogonal
\begin{equation}\label{Hermiteorthorgonal}
	(\widehat H_n(\xi),\widehat H_m(\xi))=\int_{-\infty}^{+\infty}\widehat H_n(\xi)\widehat H_m(\xi)dx=\delta_{mn},
\end{equation}
where $\delta_{mn}$ is Kronecker symbol.

We discretize the computational time interval $[0, T]$ into equally-spaced intervals $I_k:=[k\Delta t, (k+1)\Delta t]$ for $k=0, 1, \cdots, N$, where $\Delta t=T/N$.  Then, the Hermite spectral method together with backward Euler time discretization is to find approximation
\begin{equation}\label{Hermiteapprox}
	\tilde\eta^{\pm}_{m}(\xi)=\sum_{k=0}^{p}\tilde\eta_{mk}^{\pm}\widehat H_k(\lambda\xi),
\end{equation}
for $\hat\eta^{\pm}(\xi, t)$ at time $t_m=m\Delta t$ s.t.,
\begin{equation}\label{hermitespectralscheme}
	\Big(\frac{\tilde\eta^{\pm}_{m}(\xi)-\tilde\eta^{\pm}_{m-1}(\xi)}{\Delta t},\widehat H_n(\lambda\xi)\Big)=-a(\tilde\eta^{\pm}_{m}(\xi),\widehat H_n(\lambda\xi)),
\end{equation}
for all $n=0, 1, \cdots, p$. Here,
$\lambda$ is a scaling parameter to achieve resolution near $\xi=0$, and the bilinear form $a(\cdot, \cdot)$ is defined as
\begin{equation}\label{bilinearform} a(\phi(\xi),\psi(\xi))=\Big(A_s^{\mp}(\xi)\frac{d\phi(\xi)}{d\xi}, \frac{d\psi(\xi)}{d\xi}\Big)-(B_s^{\pm}(\xi)\phi(\xi). \psi(\xi)).
\end{equation}
Next, with the unknown vector denoted by $\bs U^{\pm}_m=(\tilde\eta^{\pm}_{m0}, \tilde\eta^{\pm}_{m1}, \cdots, \tilde\eta^{\pm}_{mp})^{\rm T}$, the numerical scheme \eqref{hermitespectralscheme} gives a linear system
\begin{equation}\label{matvecform}
	\mathbb D\frac{\bs U^{\pm}_m-\bs U^{\pm}_{m-1}}{\Delta t}=(\mathbb K^{\mp}+\mathbb M^{\pm})\bs U^{\pm}_m,
\end{equation}
where $\mathbb D=(D_{nk})$, $\mathbb K^{\pm}=(K_{nk}^{\pm})$, $\mathbb M=(M_{nk}^{\pm})$ are matrices with entries given by
\begin{equation}
	\begin{split}
		D_{nk}&=(\widehat H_k(\lambda \xi), \widehat H_n(\lambda \xi))=\frac{1}{\lambda}\delta_{nk},\quad K_{nk}^{\pm}=-\lambda^2\Big(A_s^{\pm}(\xi)\widehat H'_k(\lambda\xi), \widehat H'_n(\lambda\xi)	\Big),\\
		M_{nk}^{\pm}&=-(B^{\pm}_s(\xi)\widehat  H_k(\lambda\xi), \widehat H_n(\lambda\xi)).
	\end{split}
\end{equation}

By using the recurrence formula of the Hermite functions, formulations for the matrices $\mathbb K^{\pm}$, $\mathbb M^{\pm}$ can be derived analytically (see the appendix A).
\begin{figure}[ht!]
	\center
	\subfloat[$A_0^-(\xi)$, $B_0^+(\xi)$]{\includegraphics[scale=0.21]{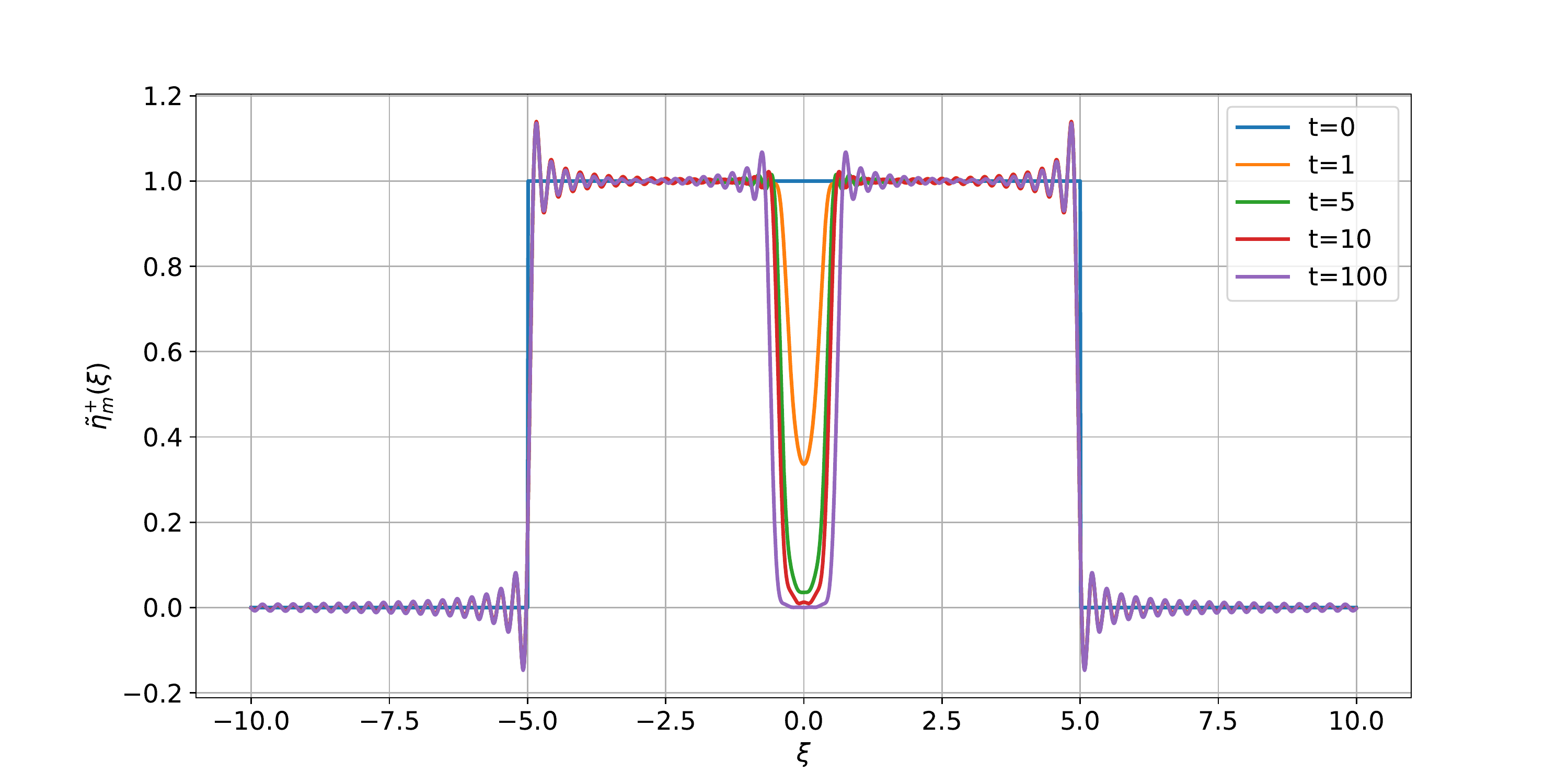}}
	\subfloat[$A_0^+(\xi)$, $B_0^-(\xi)$]{\includegraphics[scale=0.21]{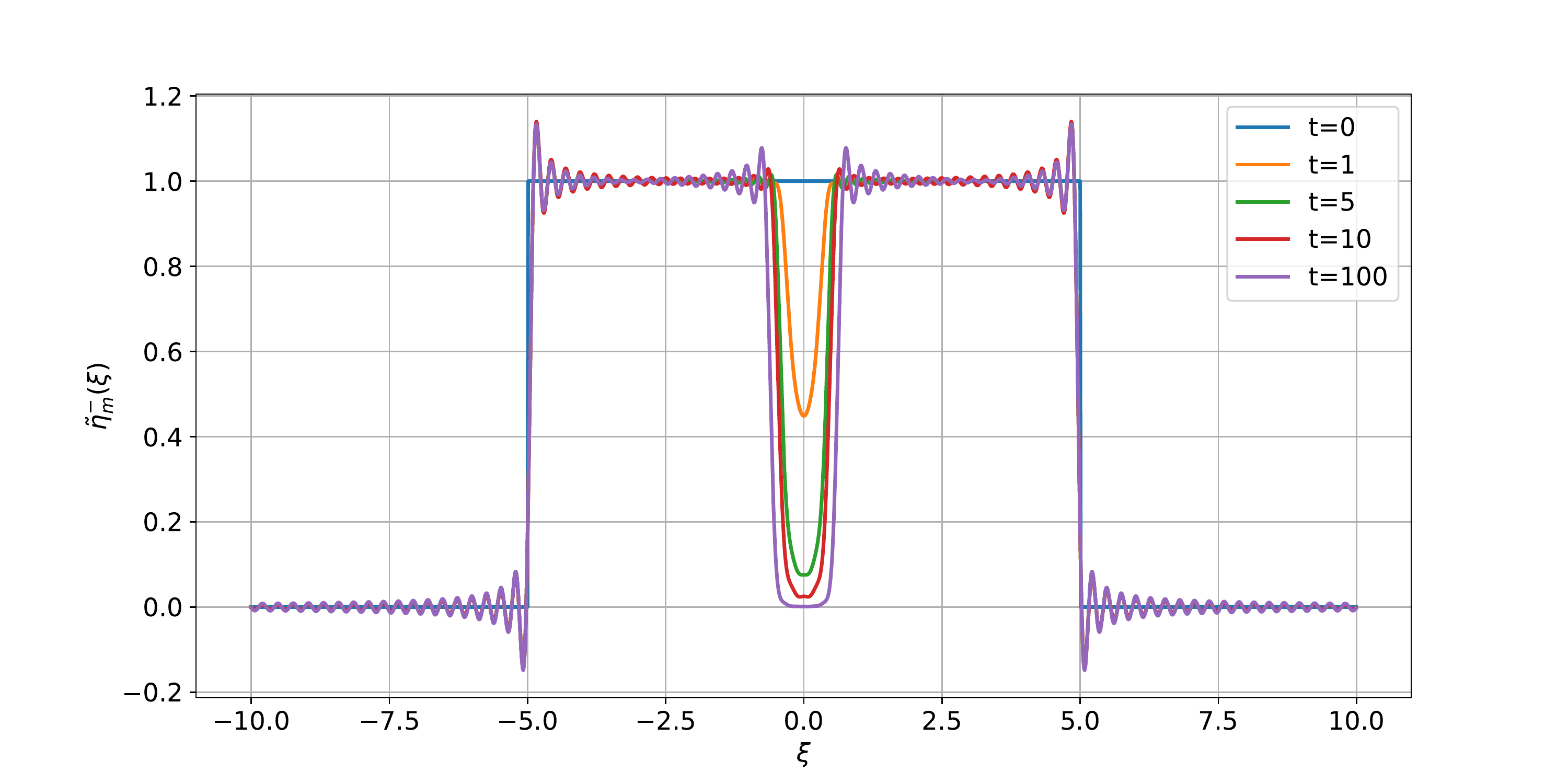}}\\
	\subfloat[$A_3^-(\xi)$, $B_3^+(\xi)$]{\includegraphics[scale=0.21]{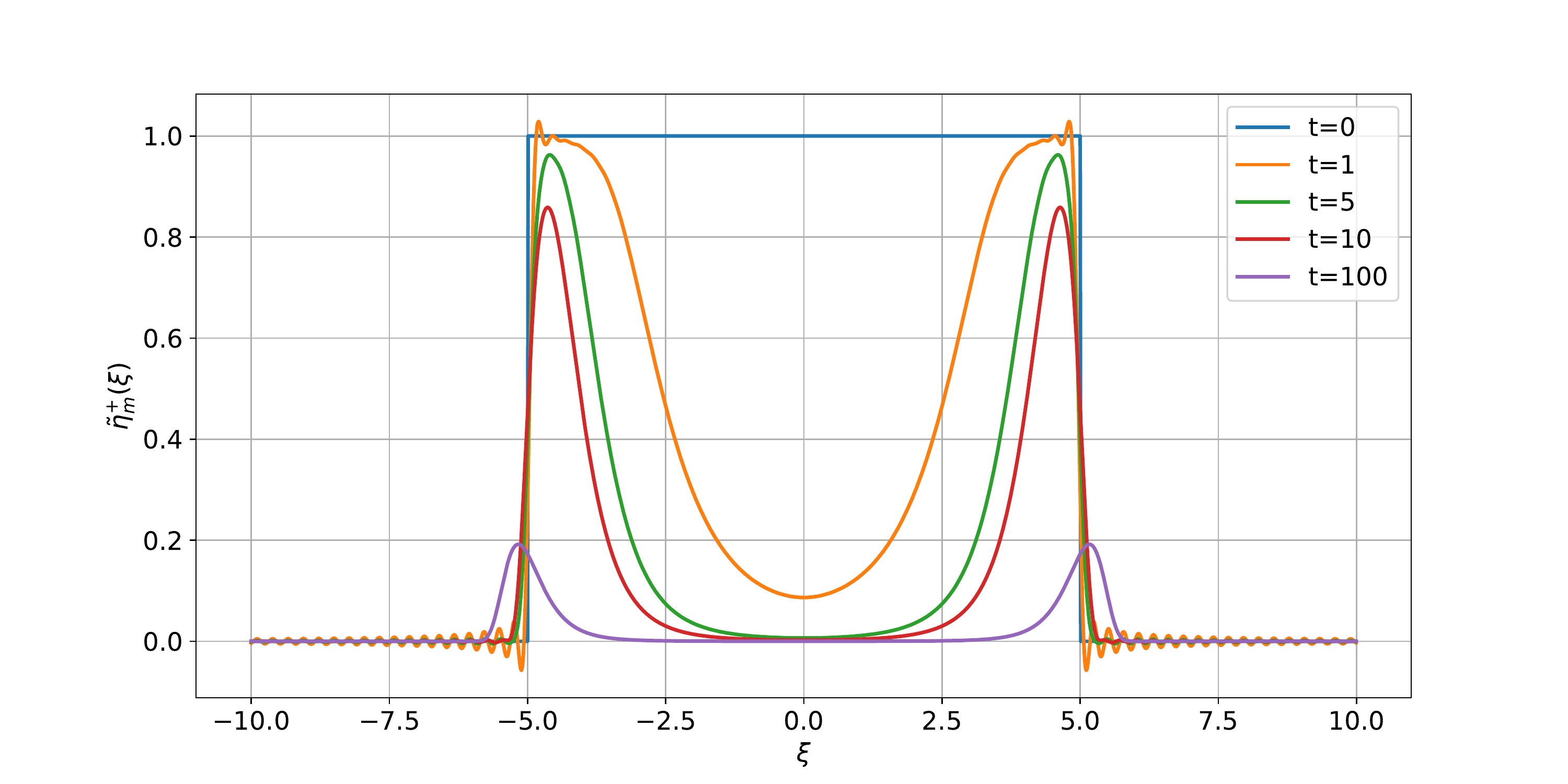}}
	\subfloat[$A_3^+(\xi)$, $B_3^-(\xi)$]{\includegraphics[scale=0.21]{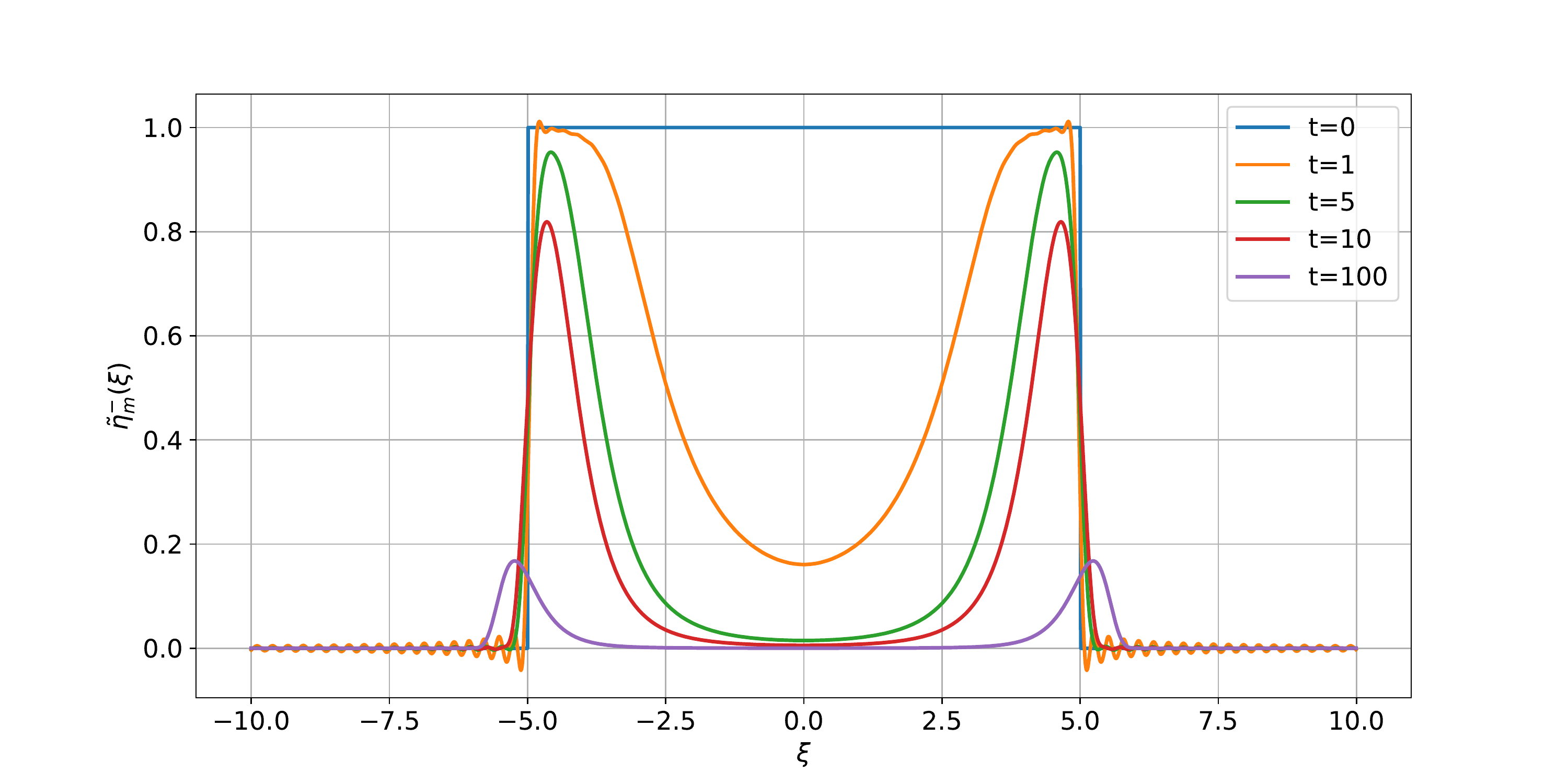}}\\
	\subfloat[$A_5^-(\xi)$, $B_5^+(\xi)$]{\includegraphics[scale=0.21]{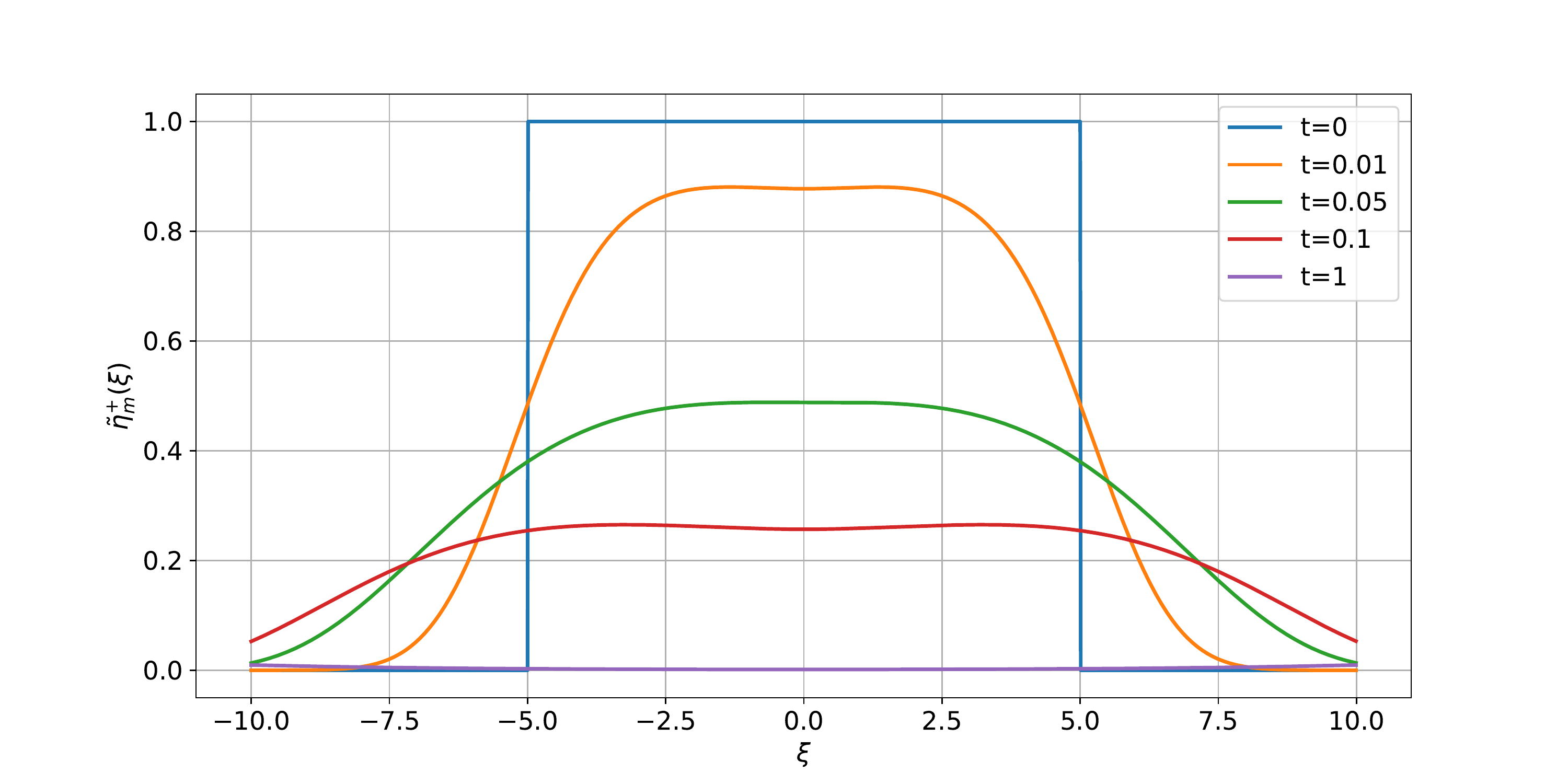}}
	\subfloat[$A_5^+(\xi)$, $B_5^-(\xi)$]{\includegraphics[scale=0.21]{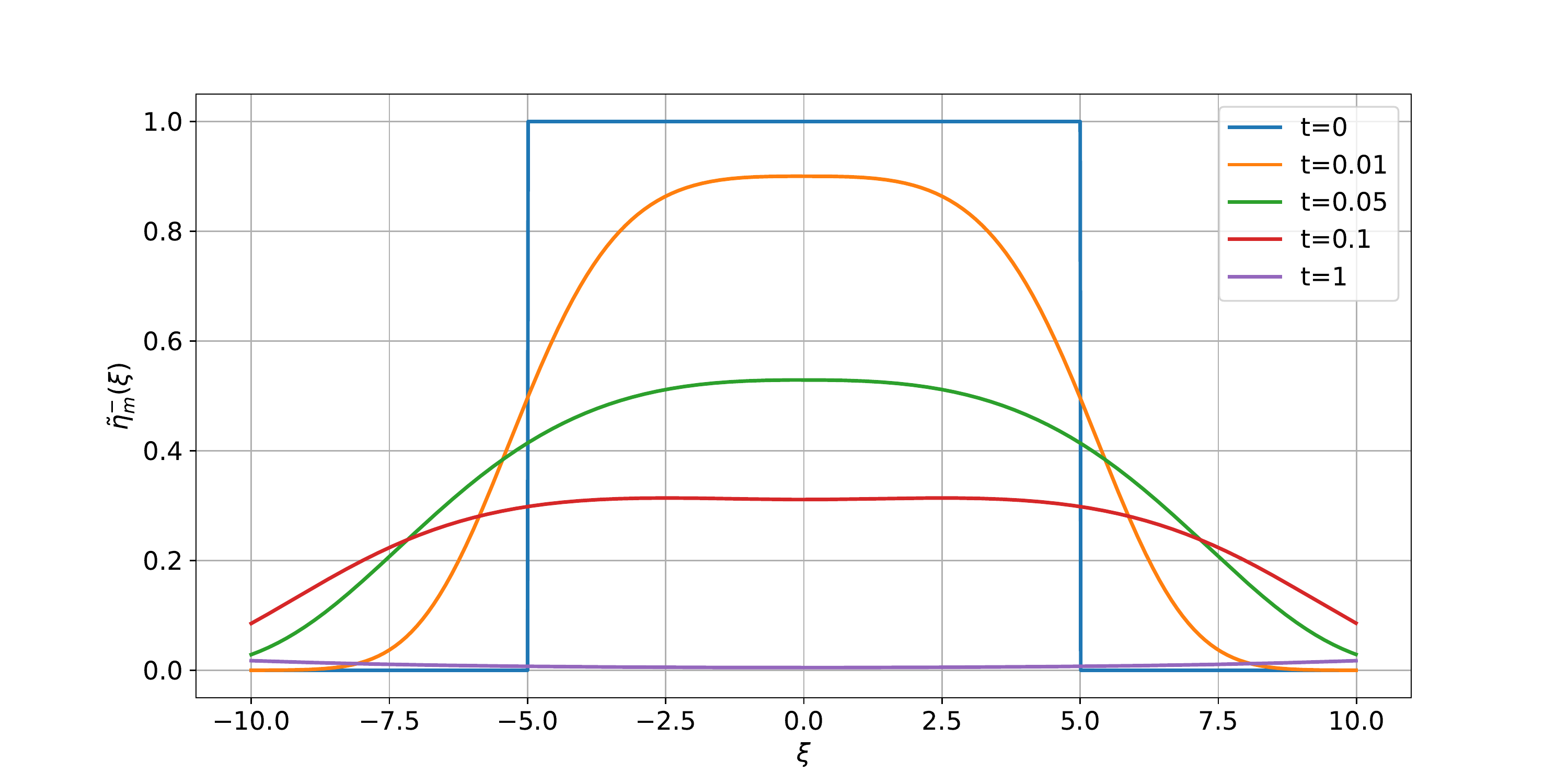}}
	\caption{Frequency domain error decay in time for a regression problem as predicted by the diffusion model \eqref{realimagdynamicsbaistrueMsDNN} for three MscaleDNNs of scale ($s=0,3,5$) with corresponding diffusion coefficients $A_s^{\pm}(\xi), B_s^{\pm}(\xi)$.}%
	\label{Ex1_1}%
\end{figure}
\subsection{Spectral bias reduction of a two layer MscaleDNN}
Some numerical examples will be presented to show the capability of the diffusion model in predicting the error dynamics of a two layer MscaleDNN. The predicted results will be compared with the training error of the two-layer MscaleDNN with a large network width and sine activation function. The spectral bias reduction phenomena of MscaleDNNs is validated by the numerical solution of the diffusion model.

{\bf Test 1 (Decaying behavior predicted by the diffusion model).} We first study the decay speed and range of the solution of the diffusion model \eqref{realimagdynamicsbaistrueMsDNN}.
Considering an initial condition for the error function in the frequency domain
\begin{equation}
	\hat\eta^{\pm}(\xi, 0)=\begin{cases}
		1, & |\xi|\leq 5,\\
		0, & |\xi|>5,
	\end{cases}
\end{equation}
we will test the diffusion model \eqref{realimagdynamicsbaistrueMsDNN} with three sets of coefficients $\{A^{\mp}_s(\xi), B_s^{\pm}(\xi)\}$, $s=0, 3, 6$. For the numerical discretization of the PDE, we take $p=100$, $\Delta t=1.0e-3$ in \eqref{Hermiteapprox}. The numerical solutions at different time $t$ are plotted in Fig. \ref{Ex1_1}. The numerical results clearly show that the initial error function decays faster over wider frequency ranges with an increasing of $s$. It is worthy to emphasize that diffusion coefficients $\{A^{\mp}_0(\xi), B_0^{\pm}(\xi)\}$ only produce fast decay in only a small neighborhood of the zero frequency, which corresponds to exactly the spectral bias of a fully connected DNN \cite{rahaman2019spectral, xu2020frequency}. These observations are consistent with the performance of the MscaleDNN, which has faster convergence  in the approximation of highly oscillated functions.

\medskip

{\bf Test 2 (Validation of error diffusion model with real MscaleDNN training).} In this test, we will show that the error dynamics of a finite but wide enough 2-layered multi-scale neural network can be predicted by the diffusion equation model quite well.

We consider a fitting problem with an objective function
\begin{equation}\label{objfun}
	f(x)=\sin a\pi x+\cos b\pi x,
\end{equation}
on the interval $[-\beta, \beta]$.
The Fourier transform of $f(x)$ with zero extension outside $[-\beta,\beta]$ is
\begin{equation*}
	\hat f(\xi)=\frac{\sin[(b+2\xi)\beta\pi]}{(b+2\xi)\pi}+\frac{\sin[(b-2\xi)\beta\pi]}{(b-2\xi)\pi}+\ri\left[\frac{\sin[(a+2\xi)\beta\pi]}{(a+2\xi)\pi}-\frac{\sin[(a-2\xi)\beta\pi]}{(a-2\xi)\pi}\right].
\end{equation*}
For the two layers multi-scale neural network, the Fourier transform of $\mathcal N_s(x,\theta)$ with zero extension outside $[-\beta,\beta]$ can be calculated as
\[
\widehat{\mathcal N}_s(\xi,\theta)=\frac{1}{\sqrt{N}}\sum\limits_{j=0}^{s}\alpha_j\sum\limits_{k=1}^{q} S_{j,k}(\xi,\theta) +\frac{1}{\sqrt{N}}\sum\limits_{j=0}^{s}\alpha_j\sum\limits_{k=1}^{q} C_{j,k}(\xi,\theta),
\]
where
\[
S_{j,k}(\xi,\theta) =\frac{-2\pi\ri\xi(e^{2\pi\ri\beta\xi}\sin(\alpha_j\theta_{jq+k}\beta-b_{jq+k})+e^{-2\pi\ri\beta\xi}\sin(\alpha_j\theta_{jq+k}\beta+b_{jq+k}))}{\alpha^2_j\theta_{jq+k}^2-4\pi^2\xi^2},
\]
and
\[
C_{j,k}(\xi,\theta)=\frac{\alpha_j\theta_{jq+k}(e^{2\pi\ri\beta\xi}\cos(\alpha_j\theta_{jq+k}\beta-b_{jq+k})-e^{-2\pi\ri\beta\xi}\cos(\alpha_j\theta_{jq+k}\beta+b_{jq+k})) }{\alpha^2_j\theta_{jq+k}^2-4\pi^2\xi^2}.
\]
\begin{figure}[ht!]
	\center
	\subfloat{\includegraphics[scale=0.21]{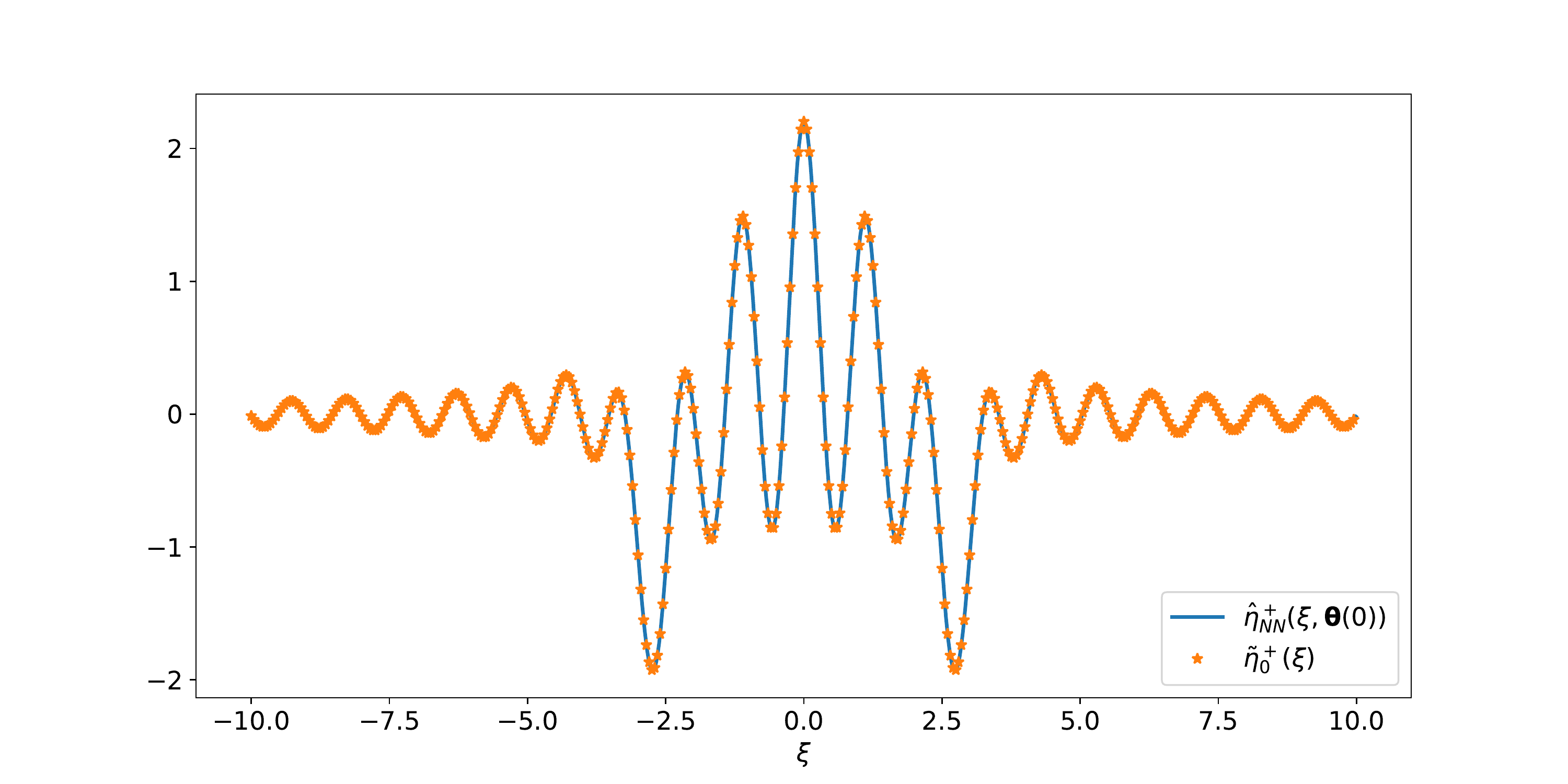}}
	\subfloat{\includegraphics[scale=0.21]{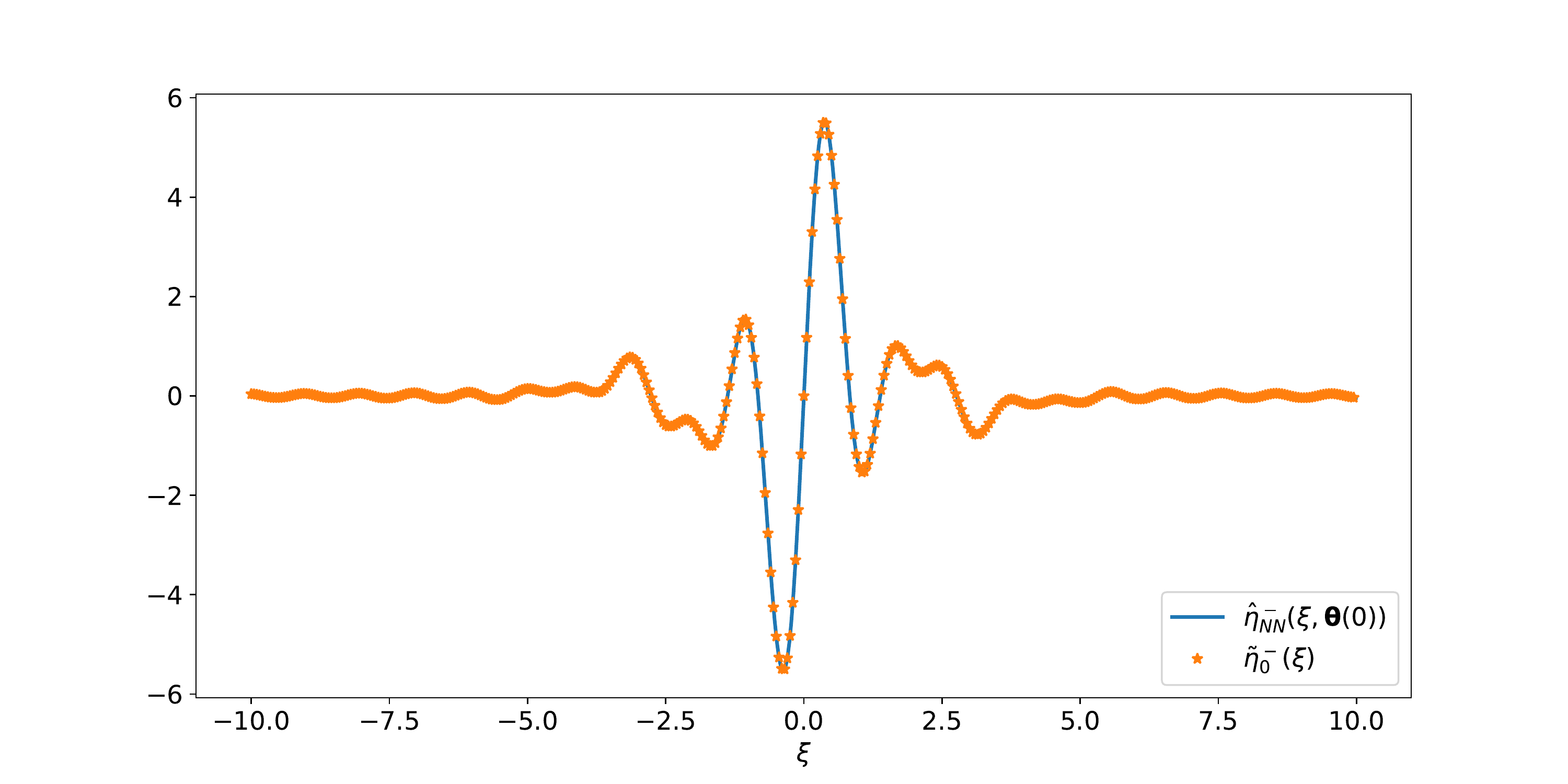}}\\
	\subfloat{\includegraphics[scale=0.21]{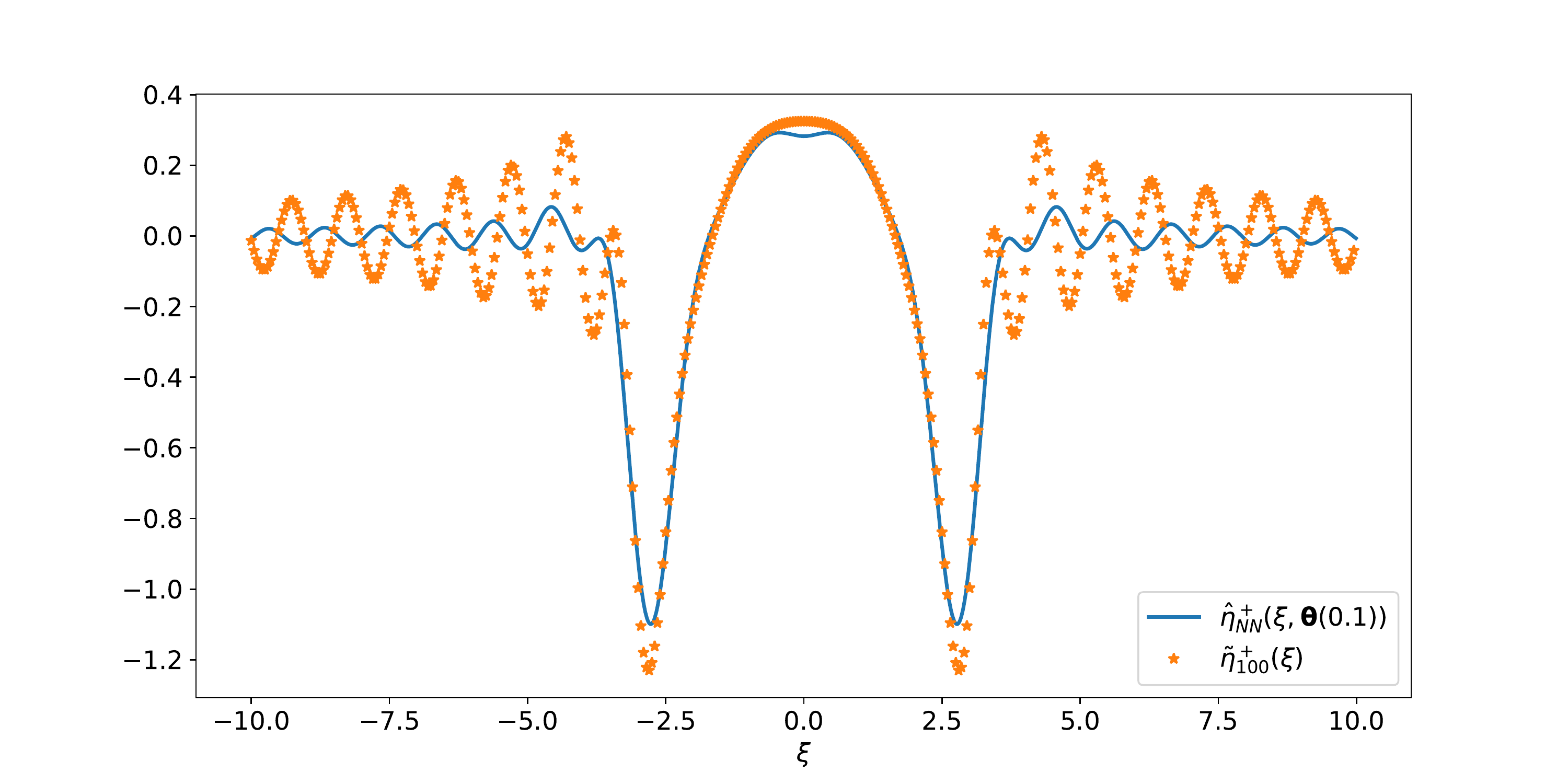}}
	\subfloat{\includegraphics[scale=0.21]{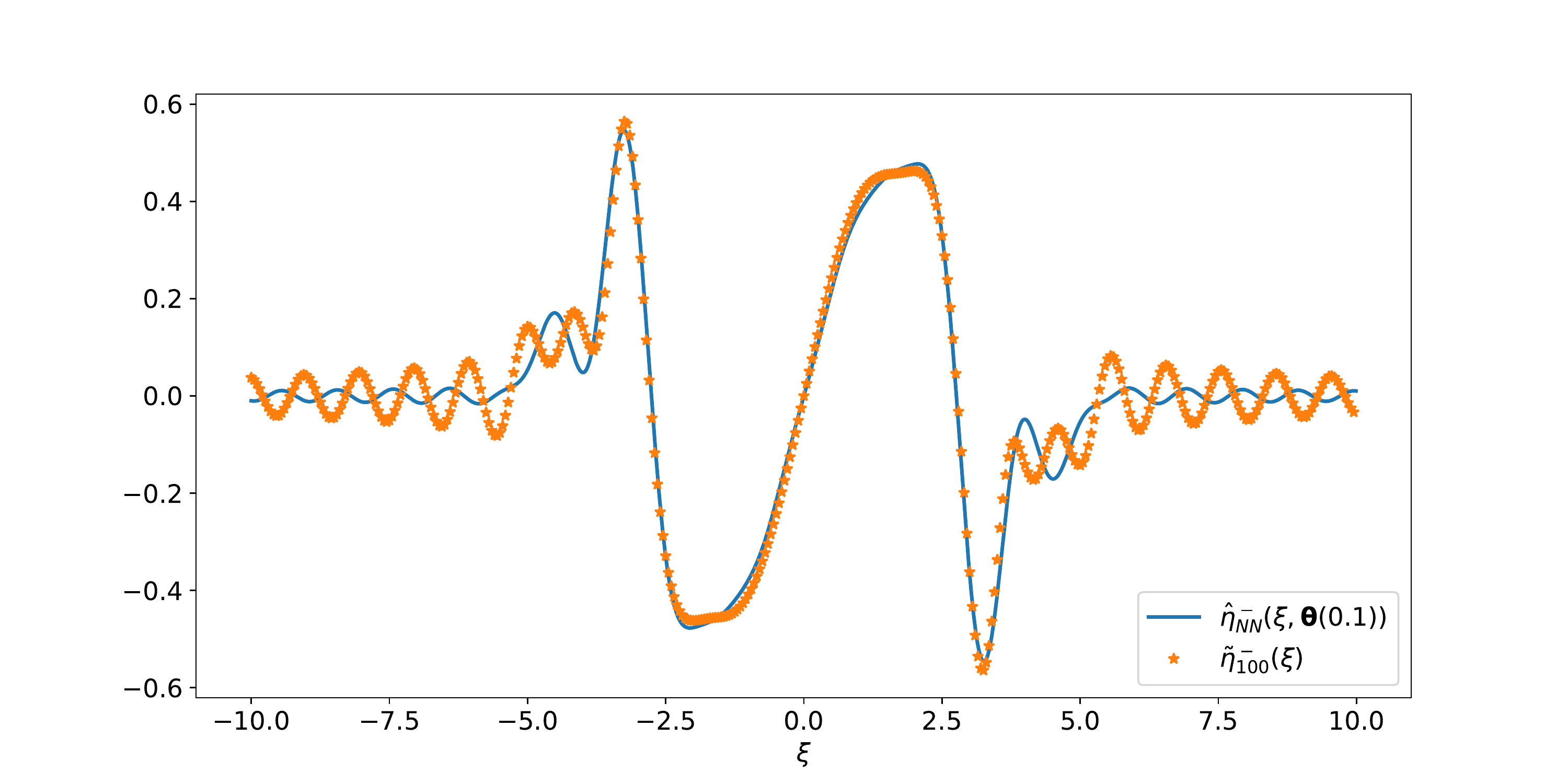}}\\
	\subfloat{\includegraphics[scale=0.21]{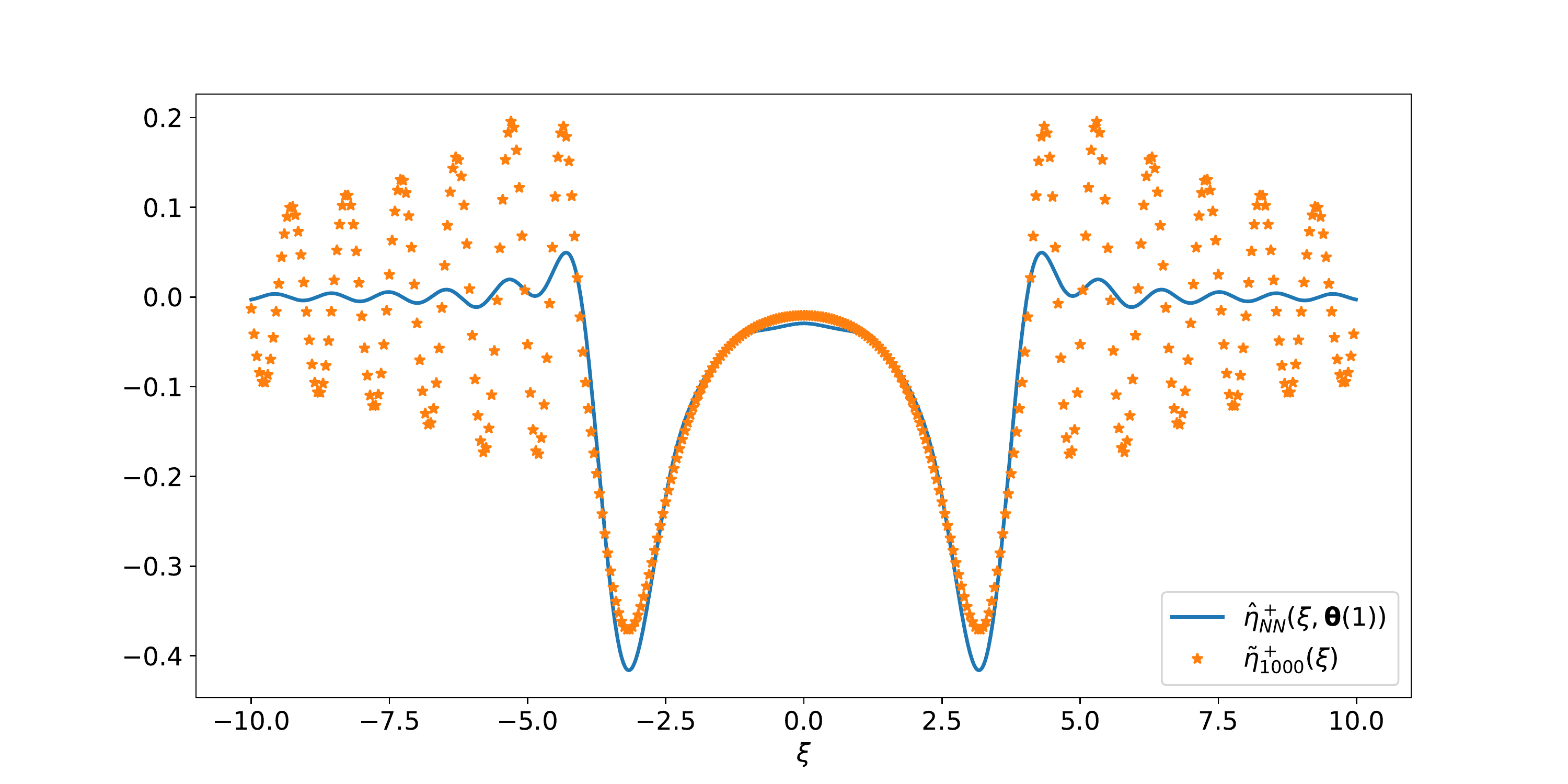}}
	\subfloat{\includegraphics[scale=0.21]{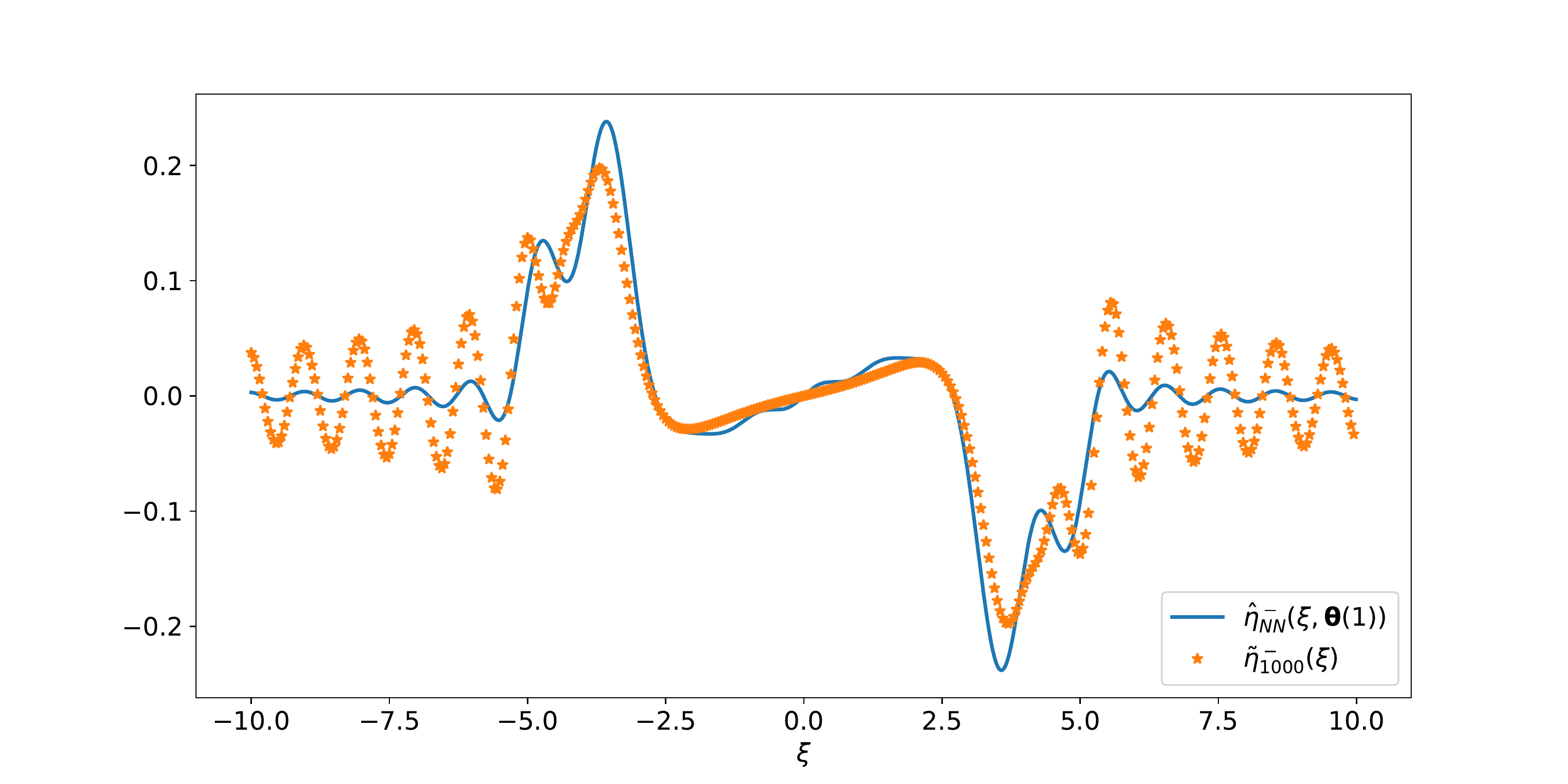}}\\
	\subfloat{\includegraphics[scale=0.21]{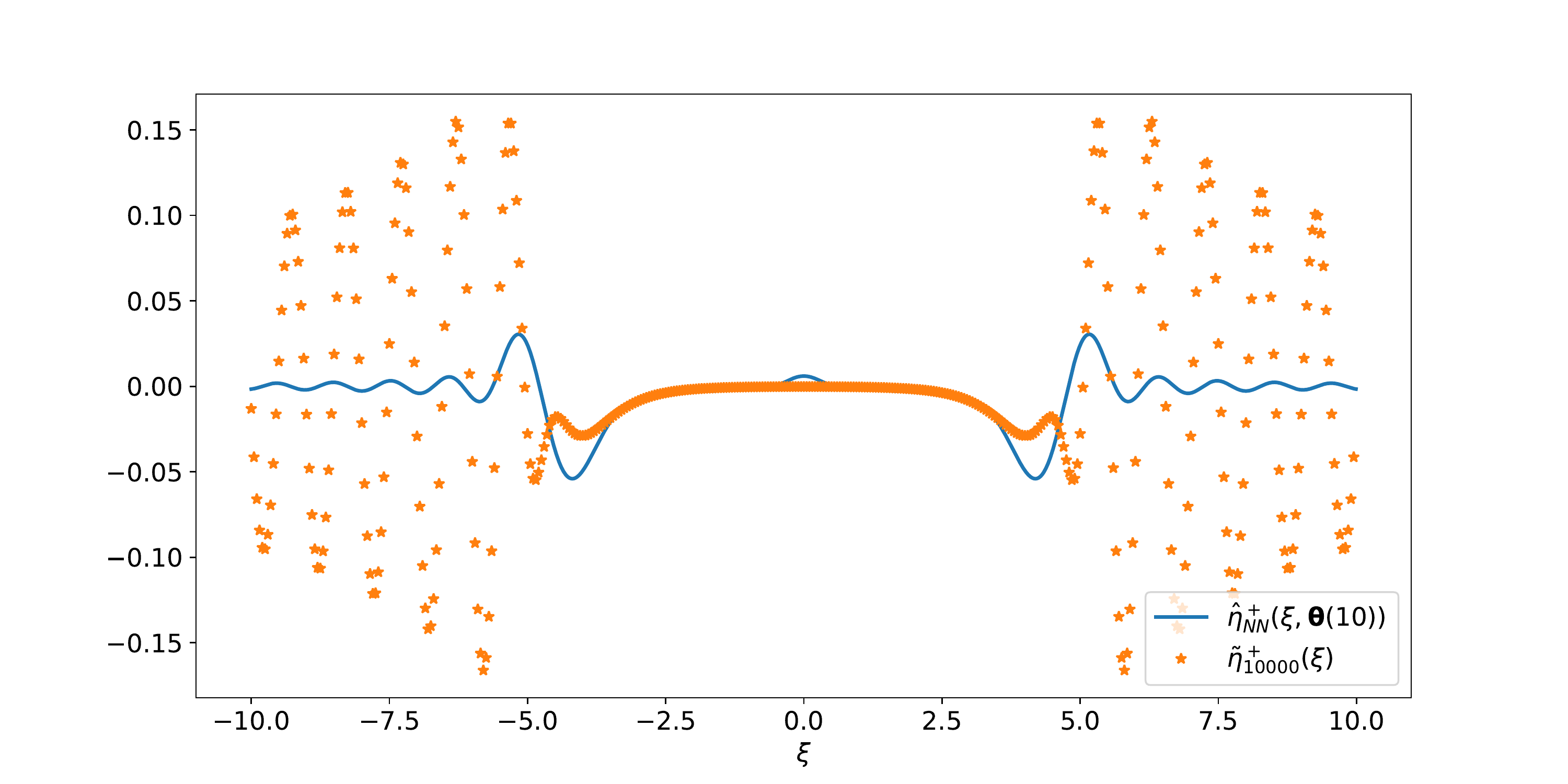}}
	\subfloat{\includegraphics[scale=0.21]{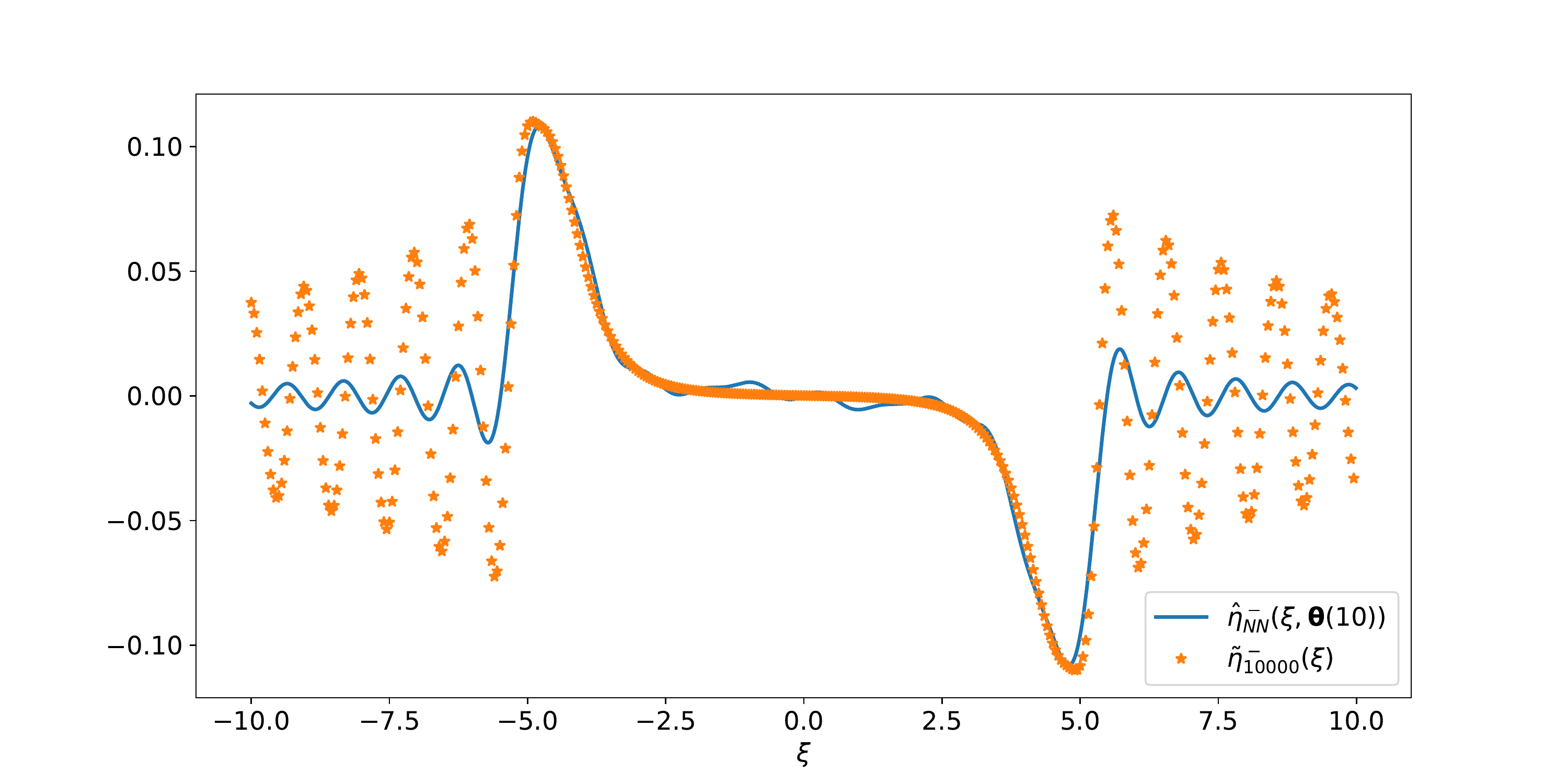}}
	\caption{Frequency domain error evolution (left - real part, right - imaginary part) in time of a 3-scale MscaleDNN with a network width $N=12,000$ (line) vs prediction by diffusion model \eqref{realimagdynamicsbaistrueMsDNN} (symbol).}%
	\label{Ex2_1}%
\end{figure}
\begin{figure}[ht!]
	\center
	\subfloat{\includegraphics[scale=0.21]{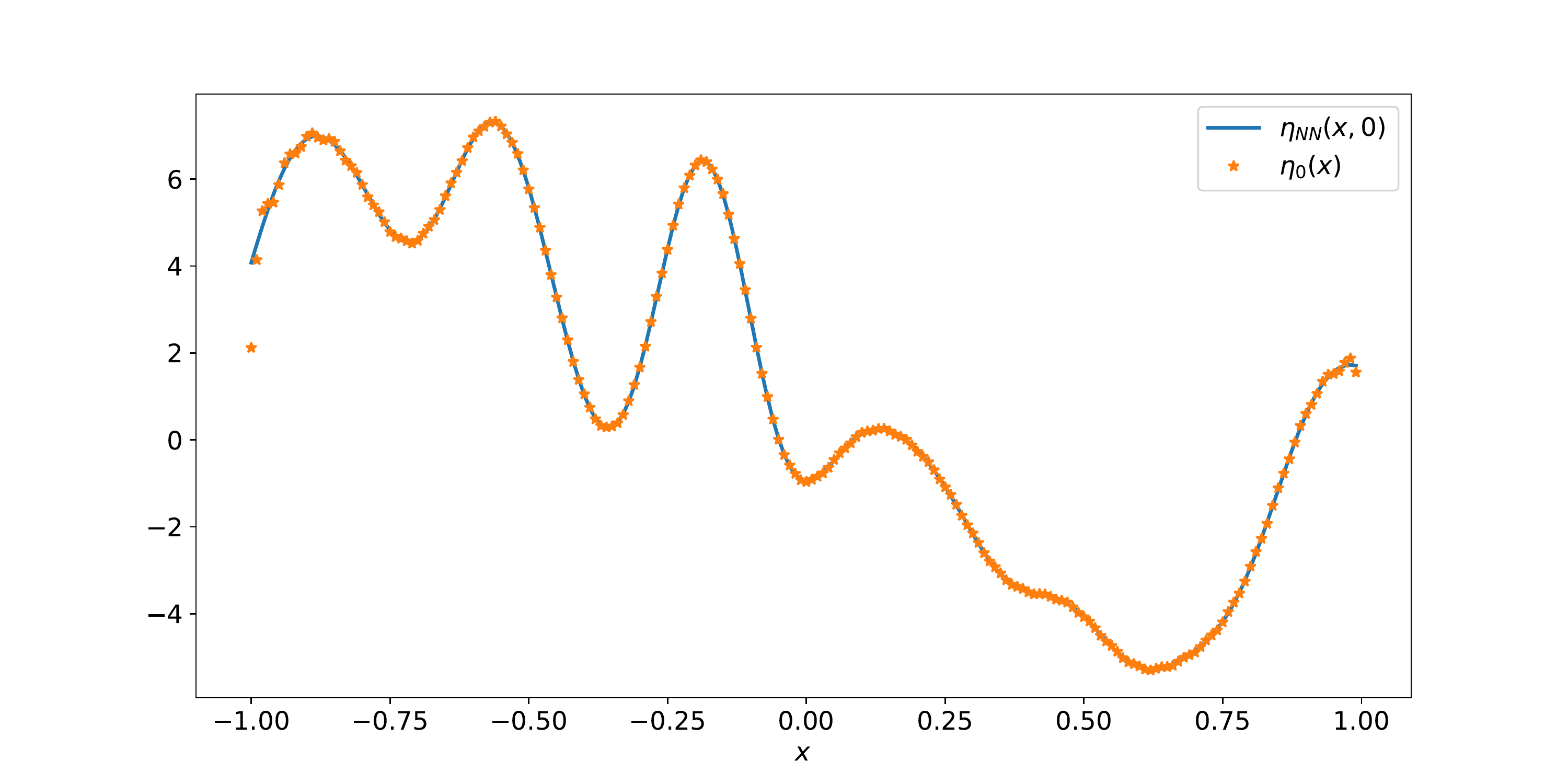}}
	\subfloat{\includegraphics[scale=0.21]{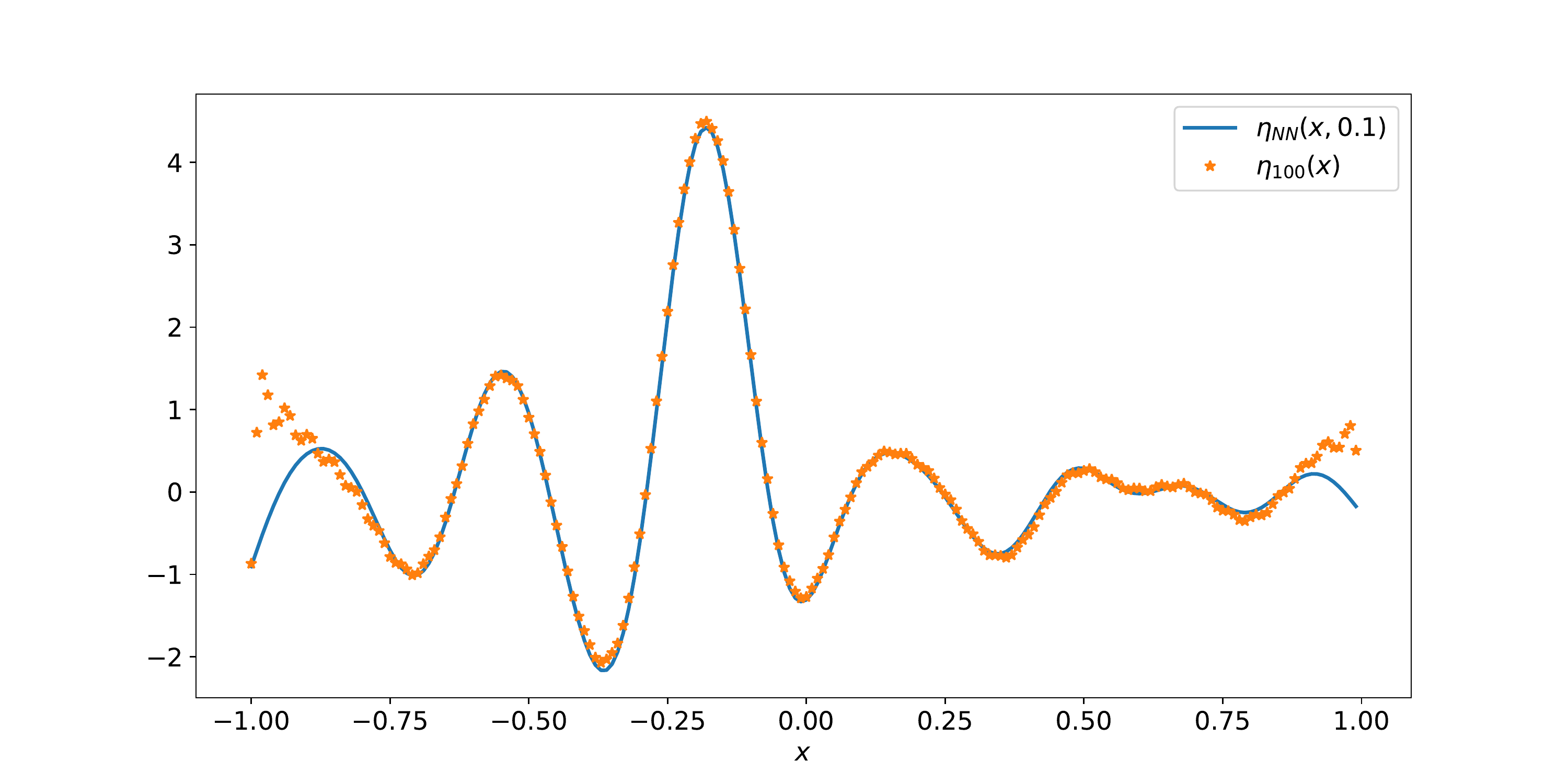}}\\
	\subfloat{\includegraphics[scale=0.21]{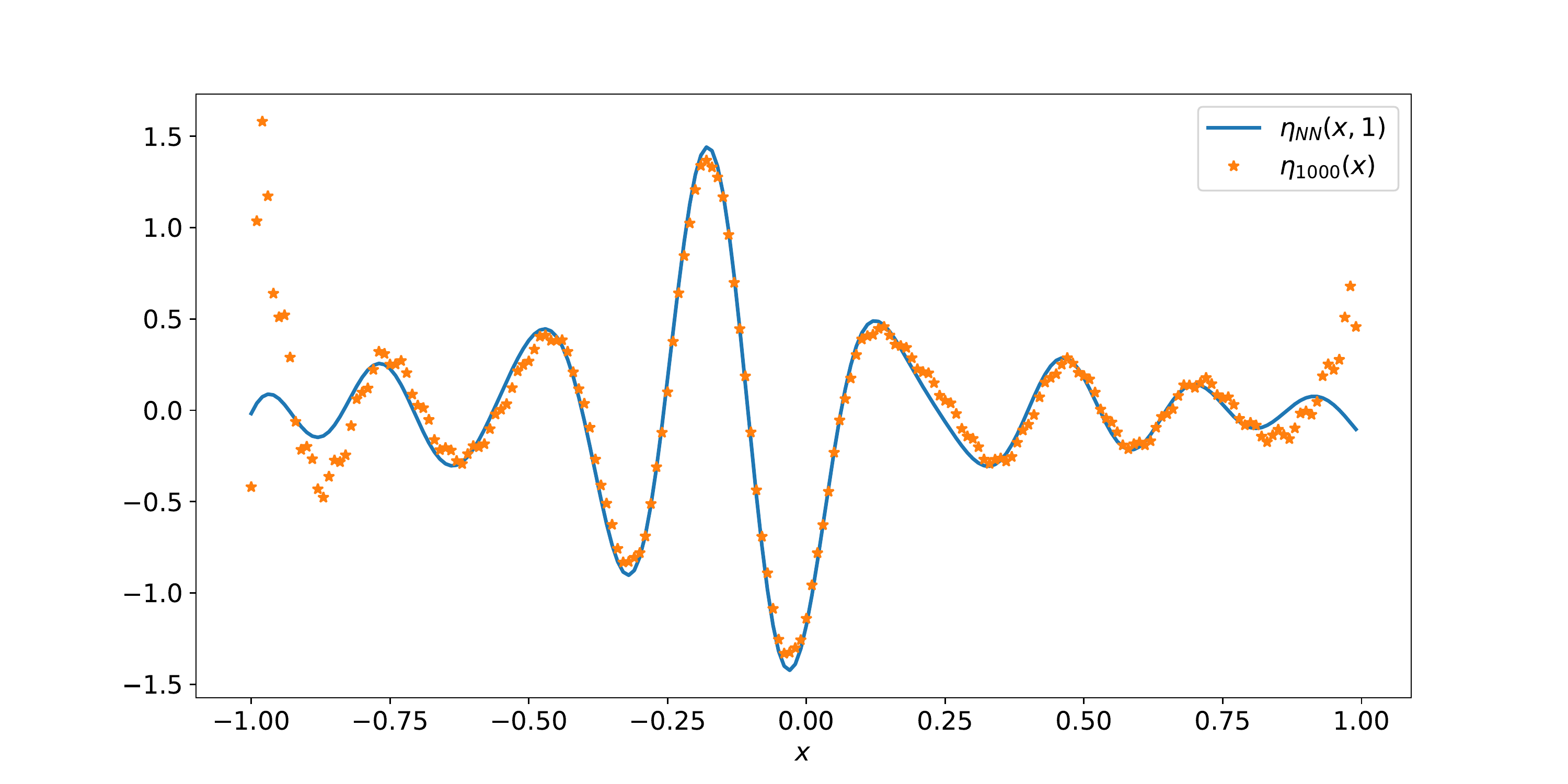}}
	\subfloat{\includegraphics[scale=0.21]{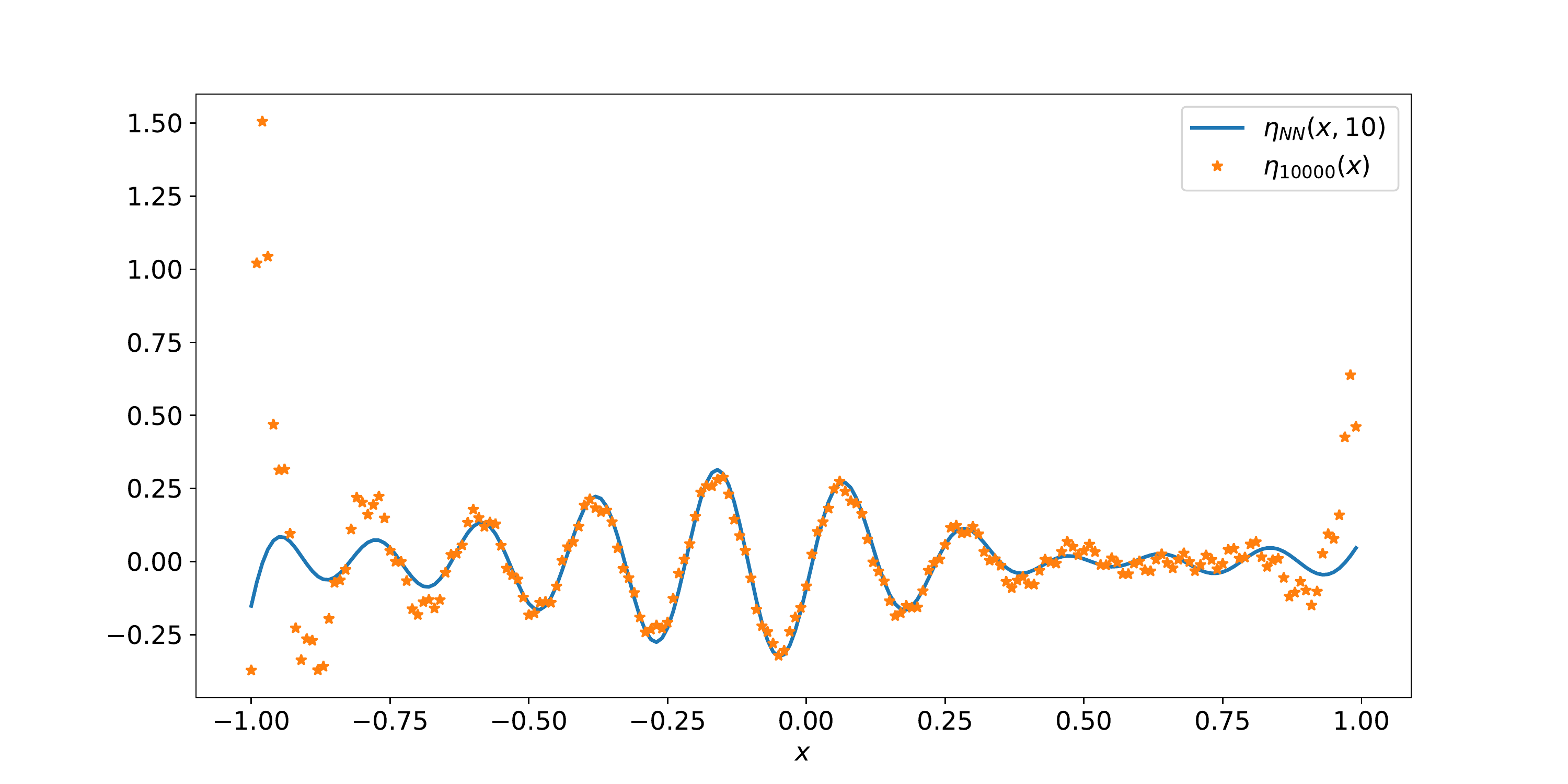}}
	\caption{Physical domain error evolution in time of a 3-scale MscaleDNN with a network width $N=12,000$ (line) vs prediction by diffusion model \eqref{realimagdynamicsbaistrueMsDNN} (symbol).}%
	\label{Ex2_2}%
\end{figure}
We will show that the error $\hat\eta_{NN}(\xi, \theta)=\widehat{\mathcal N}_s(\xi,\theta)-\hat f(\xi)$ of the MscaleDNN by the gradient descent learning agrees with that predicted by the diffusion equation \eqref{realimagdynamicsbaistrueMsDNN}. We take $a=4.2$, $b=5.8$, $\beta=1$ and the initial errors are given by $\eta_{NN}(x, \theta_0)={\mathcal N}_s(x,\theta_0)-f(x)$  with parameters initialized by sampling from independent random variables of normal distribution. In the gradient descent training for the $\mathcal N_s(x,\theta)$, the training data set consists of $2000$ uniformly distributed points in $[-\beta, \beta]$ and learning rate $\tau=1.0e-3$ is adopted. In this example, a two layers neural network with  $m=12,000$, $\alpha_j=2^j$ and scale $s=3$ is tested and the training is performed in full batch.

Meanwhile, in the Fourier spectral domain, the diffusion equation \eqref{realimagdynamicsbaistrueMsDNN} with initial function $\hat\eta(\xi, \theta_0)=\widehat{\mathcal N}_s(\xi,\theta_0)-\hat f(\xi)$ will be solved with a $p$-th order the Hermite spectral method introduced above.  We take $p=300$ and $\Delta t=\tau$ in the discretization.

The Fourier transform of $\eta_{NN}(x,\theta(t))$, denoted by
$$\hat\eta_{NN}(x,\theta(t))=\hat \eta_{NN}^+(\xi,\theta(t))+\ri \eta_{NN}^-(\xi,\theta(t))$$
are compared with $\tilde \eta^{\pm}_m(\xi)$ at $t=m\Delta t$, see Fig. \ref{Ex2_1}. Although many approximations have been used in deriving the diffusion model, the results show that the prediction produced by the diffusion model captured the main features of the error over a long time training process.

On the other hand, we can also compare the training error with the diffusion model prediction in the physical domain. Using the fact that \cite{gradshteyn2014table,li2022efficient}
\begin{equation}
	\mathcal F^{-1}[\widehat H_k(\xi)](x)=\int_{-\infty}^{+\infty}\widehat H_k(\xi)e^{2\ri\pi\xi x}{\rm d}\xi=\sqrt{2\pi}\ri^k\widehat H_k(2\pi\xi),
\end{equation}
the Hermite approximation of the error predicted by the diffusion model, i.e.,
\begin{equation}
	\tilde\eta_m^{\pm}(\xi)=\sum_{k=0}^{p}\tilde\eta_{mk}\widehat H_k(\lambda\xi),
\end{equation}
can be analytically transformed back to the physical domain as
\begin{equation}
	\eta_m^{\pm}(x):=\mathcal F^{-1}[\tilde\eta^{\pm}_m](x)=\sum_{k=0}^{p}\tilde\eta_{mk}\int_{-\infty}^{+\infty}\widehat H_k(\lambda\xi)e^{2\ri\pi x\xi}{\rm d}\xi=\frac{\sqrt{2\pi}}{\lambda}\sum_{k=0}^{p}\tilde\eta_{mk}\ri^k\widehat H_k\Big(\frac{2\pi x}{\lambda}\Big).
\end{equation}
Then, in the physical domain the errors $\eta_{NN}(x,\theta(t_m))$ from the MscaleDNN training and $\eta_m(x)=\eta^+_m(x)+\ri \eta^-_m(x)$ predicted by the diffusion equation can be  compared in Fig. \ref{Ex2_2}. Clearly, the evolution of the errors matches quite well in physical domain. It is worthy to point out that the fitting domain $\Omega=[-1, 1]$ is not large. However, the diffusion model can still be a satisfactory predictor for the real error through the training of the MScaleDNN with a large enough network width.

\medskip

{\bf Test 3 (Reduction of spectral bias predicted by error diffusion model).}
With the confirmation of predicting capability of the diffusion equation model \eqref{realimagdynamicsbaistrueMsDNN} for the error decay of the MscaleDNN with a large enough network width, we will use the model to demonstrate the spectral bias reduction of MscaleDNNs with increasing scales.

Again, We set the network width at $m=12000$, and $a=4.2$, $b=5.8$ and the initial errors are given by $\hat\eta(\xi, \theta_0)=\widehat{\mathcal N}_s(\xi,\theta_0)-\hat f(\xi)$ with parameters initialized by sampling from independent random variables of normal distribution. In the Hermite spectral method approximation of the diffusion equation \eqref{realimagdynamicsbaistrueMsDNN}, we take $p=300$ and $\Delta t=1.0e-3$. The numerical solution of the diffusion equations at different time $t$ for a standard fully connected network (FCN) corresponding to coefficients $\{A_0^{\pm}(\xi), B_0^{\mp}\}$ and a 3-scales MscaleDNN corresponding to coefficients $\{A_3^{\pm}(\xi), B_3^{\mp}\}$  are plotted in Fig. \ref{Ex2_1}-Fig. \ref{Ex2_2}. We can see clearly that FCN with diffusion coefficients $\{A_0^{\pm}(\xi), B_0^{\mp}\}$ only produce decay in a very small neighborhood of the zero frequency while the 3-scale MscaleDNN with coefficients $\{A_3^{\pm}(\xi), B_3^{\mp}\}$ produce much faster decay in a larger frequency interval. Although the initial errors at $t=0$ are different, $\widehat{\mathcal N}_0(\xi,\theta_0)-\hat f(\xi)$ for FCN and $\widehat{\mathcal N}_3(\xi,\theta_0)-\hat f(\xi)$ for 3-scale MscaleDNN, the numerical results all verify that multi-scale neural networks has better performance in spectral bias reduction compared with the FCN.

\begin{figure}[ht!]
	\center
	\subfloat[real part]{\includegraphics[scale=0.21]{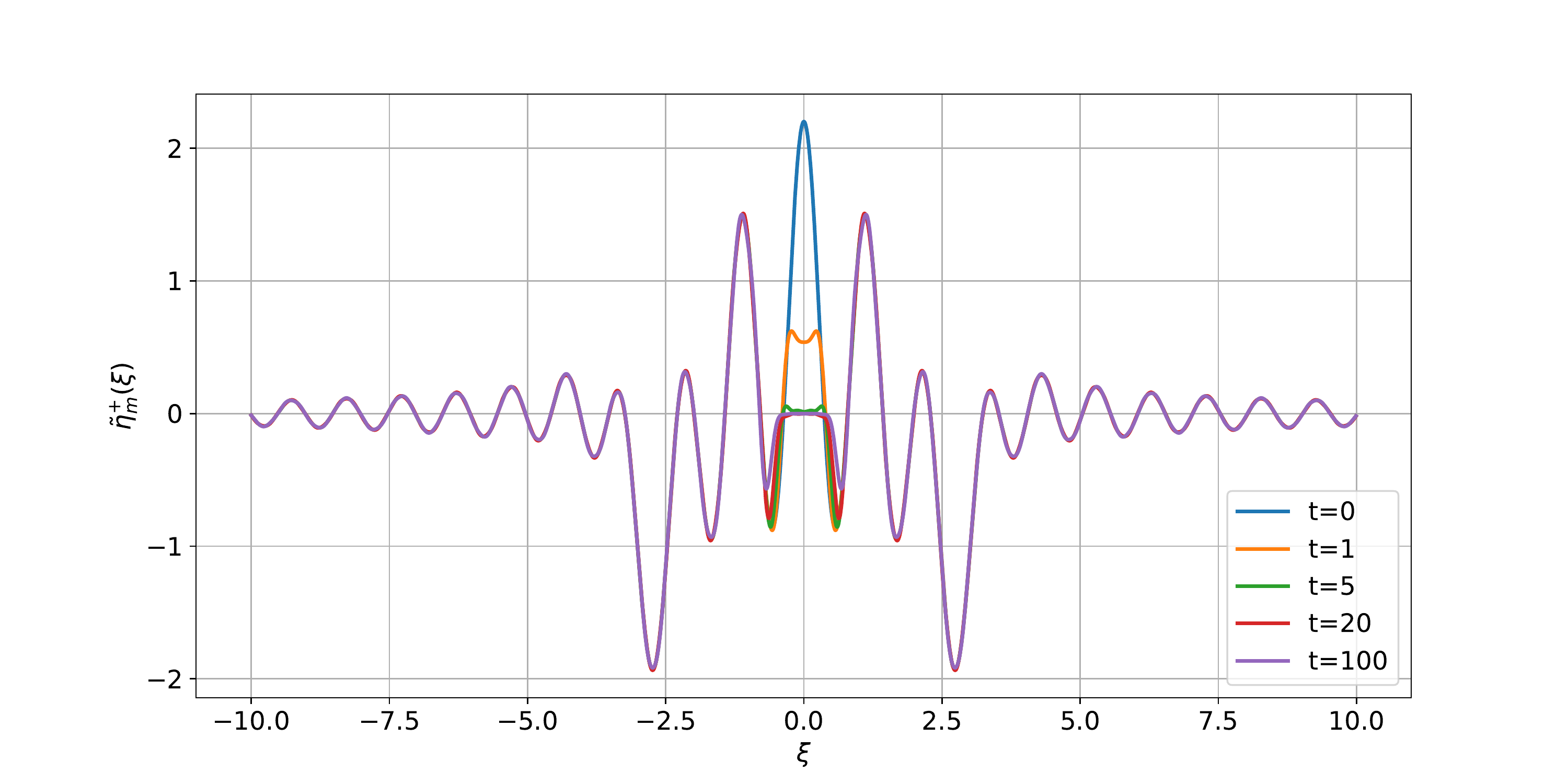}}
	\subfloat[imaginary part]{\includegraphics[scale=0.21]{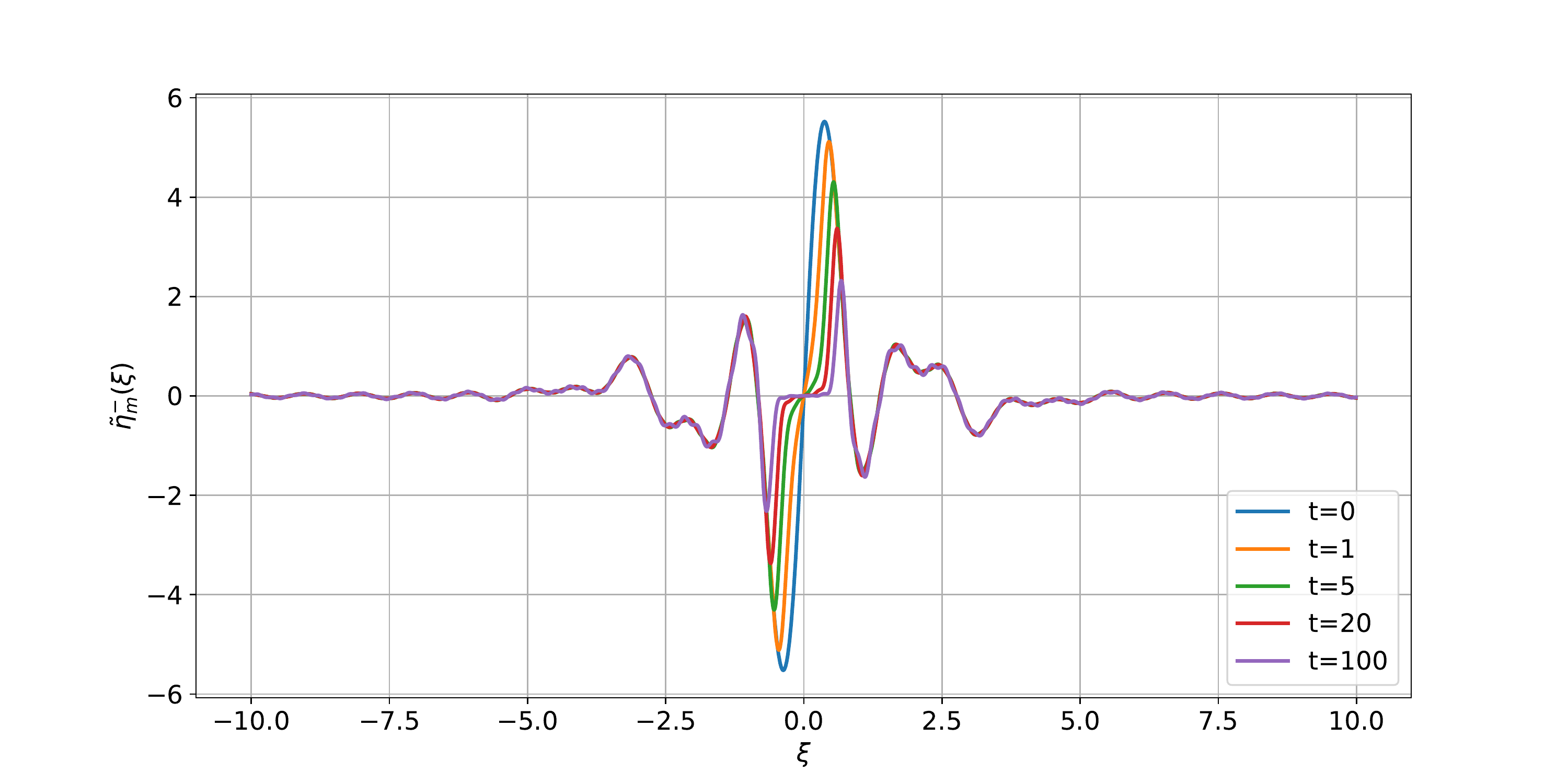}}
	\caption{Frequency domain error decay in time predicted by \eqref{realimagdynamicsbaistrueMsDNN} for a FCN corresponding to coefficients $\{A_0^{\pm}(\xi), B_0^{\mp}\}$.}%
	\label{Ex2_3}%
\end{figure}
\begin{figure}[ht!]
	\center
	\subfloat[real part]{\includegraphics[scale=0.21]{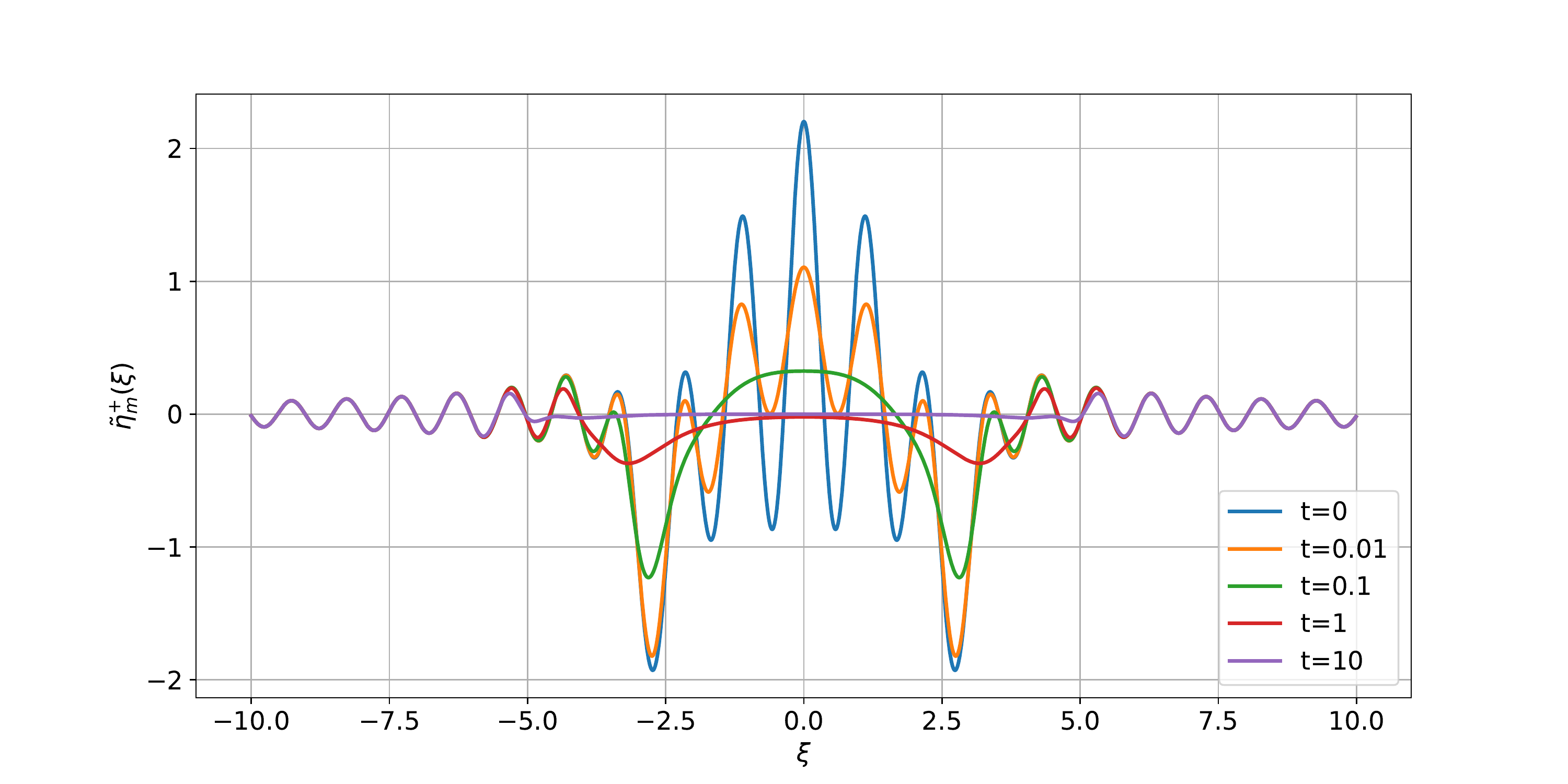}}
	\subfloat[imaginary part]{\includegraphics[scale=0.21]{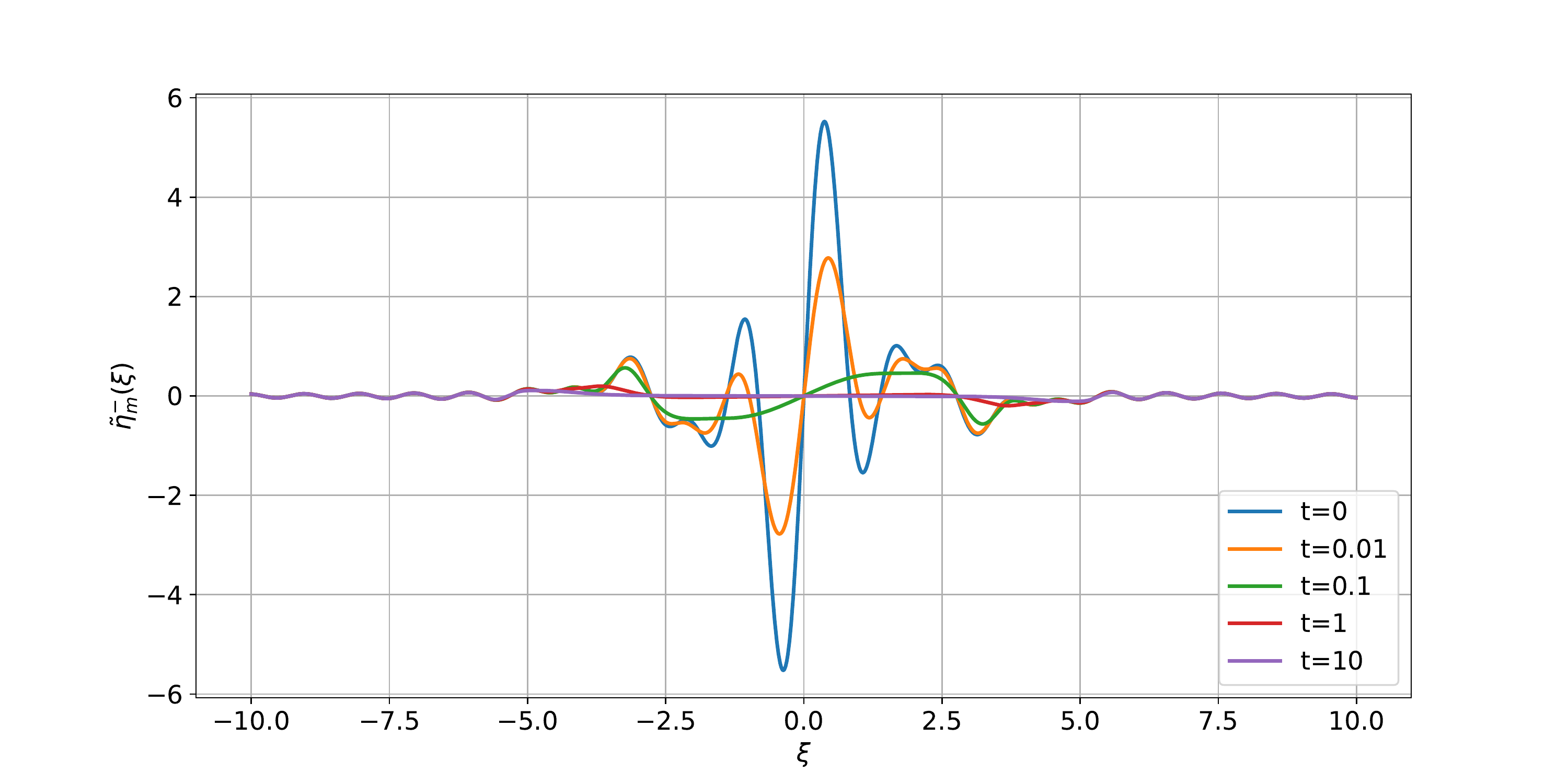}}
	\caption{Frequency domain error decay in time predicted by \eqref{realimagdynamicsbaistrueMsDNN} for a 3-scale MscaleDNN corresponding to coefficients $\{A_3^{\pm}(\xi), B_3^{\mp}\}$.}%
	\label{Ex2_4}%
\end{figure}



\section{Conclusion and future work}
In this paper, we investigated the convergence and spectral bias reduction properties of a two-layer multi-scale neural network for regression problems by deriving diffusion equation models in the frequency domain for predicting its error evolution. With the sine activation function, the gradient descent learning of MscaleDNNs leads to the diffusion equation models for the error assuming that the width of the neural network goes to infinity, the learning rate to zero and the fitting domain to the whole space. The diffusion coefficients of the diffusion equations are shown to have  wider support in the frequency domain with more scales used in the MscaleDNNs, resulting in a reduction of spectral bias for the MscaleDNNs. This is consistent with the performance of the MscaleDNN with faster convergence in approximating highly oscillated functions from various applications. Moreover, the derived diffusion equation can predict the convergence of the MscaleDNNs learning algorithm even with a finite and reasonably wide network in a finite domain.



The analysis of the MScaleDNNs with more layers, and other popular activation functions, e.g., ReLU, Sigmoid, etc, will be studied following a similar approach of this paper.

\section*{Acknowledgments}
The first and second authors acknowledges the  support provided by NSFC (grant 12022104).  W. Z. Zhang acknowledges the  support provided by NSFC (grant 12201603). The work of W. Cai is supported by the US National Science
Foundation grant DMS-2207449.

\begin{appendix}
	\section{Analytic formula for the computation of matrices $\mathbb K^{\pm},\mathbb M^{\pm}$}
	We first recall the recurrence formulas (cf. \cite{ShenTaoWang2011})
	\begin{equation}\label{H_nDS}
		\begin{split}
			\widehat H_0(x)=\pi^{-1/4}e^{-x^2/2},\quad \widehat H_1(x)=\sqrt{2}\pi^{-1/4}xe^{-x^2/2},\\
			\widehat H_{n+1}(x)=x\sqrt{\frac{2}{n+1}}\widehat H_n(x)-\sqrt{\frac{n}{n+1}}\widehat H_{n-1}(x)=0,\quad n\geq 1,
		\end{split}
	\end{equation}
	\begin{equation}\label{H_nDS1}
		\begin{split}
			\widehat H'_0(x)=-\frac{\pi^{-1/4}}{2}xe^{-x^2/2}=-\frac{\sqrt{2}}{2}\widehat H_1(x),\\
			\widehat H'_n(x)=\sqrt{2n}\widehat H_{n-1}(x)-x\widehat H_n(x)=\sqrt{\frac{n}{2}}\widehat H_{n-1}(x)-\sqrt{\frac{n+1}{2}}\widehat H_{n+1}(x),\quad n\geq 1,
		\end{split}
	\end{equation}
	of the Hermite functions $\widehat H_n(x)$.

	Then, by the recurrence formula \eqref{H_nDS1}, we have
	\begin{equation}
		\begin{split}
			&\widehat H'_0(x)\widehat H'_0(x)=\frac{1}{2}\widehat H_1(x)\widehat H_1(x), \\
			&\widehat H'_0(x)\widehat H'_n(x)=-\frac{\sqrt{n}}{2}\widehat H_1(x)\widehat H_{n-1}(x)+\frac{\sqrt{n+1}}{2}\widehat H_1(x)\widehat H_{n+1}(x),\quad n\geq 1.
		\end{split}
	\end{equation}
	and
	\begin{equation}
		\begin{aligned}
			& \widehat H'_k(x)\widehat H'_n(x)\\
   =&\left[\sqrt{\frac{k}{2}}\widehat H_{k-1}(x)-\sqrt{\frac{k+1}{2}}\widehat H_{k+1}(x)\right]\left[\sqrt{\frac{n}{2}}\widehat H_{n-1}(x)-\sqrt{\frac{n+1}{2}}\widehat H_{n+1}(x)\right]\\
			=&\frac{\sqrt{nk}}{2}\widehat H_{k-1}(x)\widehat H_{n-1}(x)-\frac{\sqrt{(n+1)k}}{2}\widehat H_{k-1}(x)\widehat H_{n+1}(x)\\
			-&\frac{\sqrt{n(k+1)}}{2}\widehat H_{k+1}(x)\widehat H_{n-1}(x)+\frac{\sqrt{(n+1)(k+1)}}{2}\widehat H_{k+1}(x)\widehat H_{n+1}(x),
		\end{aligned}
	\end{equation}
	for all $n, \;\;\; k\geq 1$. Therefore,
	\begin{equation}
		K_{00}^{\pm}=\frac{1}{2}C_{11}^{\pm},\quad K_{0n}^{\pm}=K_{n0}^{\pm}=-\frac{\sqrt{n}}{2}C_{1,n-1}^{\pm}+\frac{\sqrt{n+1}}{2}C_{1,n+1}^{\pm},\quad n\geq 1,
	\end{equation}
	where $C^{\pm}_{nk}=-\lambda^2\int_{-\infty}^{+\infty} A^{\pm}_s(\xi)\widehat H_k(\lambda\xi) \widehat H_n(\lambda\xi)d\xi.$ Otherwise, for all $n, k\geq 1$,
	\begin{equation}
		\begin{split}
			K_{nk}^{\pm}=&\frac{\sqrt{n}}{2}\big(\sqrt{k}C^{\pm}_{n-1,k-1}-\sqrt{k+1}C^{\pm}_{n-1,k+1}\big) \\
     &-\frac{\sqrt{n+1}}{2}\big(\sqrt{k}C^{\pm}_{n+1,k-1}-\sqrt{k+1}C^{\pm}_{n+1,k+1}\big).
		\end{split}
	\end{equation}
	Noting that
	\begin{equation}
		M^{\pm}_{nk}=-\int_{-\infty}^{+\infty}B_s^{\pm}(\xi)\widehat H_k(\lambda\xi)\widehat H_n(\lambda\xi)d\xi,
	\end{equation}
	and $A_s^{\pm}(\xi)$, $B_s^{\pm}(\xi)$ are linear combination of Gaussian functions as presented in \eqref{DNNcoefficients}, the computation of $C^{\pm}_{nk}$ and $M^{\pm}_{nk}$ can be reduced to compute the weighted inner products
\begin{equation}\label{generalintegral}
     \begin{split}
		I_{nk}(\tau) & =\int_{-\infty}^{+\infty}\widehat H_n(x)\widehat H_k(x)e^{-\tau x^2}dx \\
     &=\frac{1}{\sqrt{\tau+1}}\int_{-\infty}^{+\infty}\widetilde H_n\Big(\frac{y}{\sqrt{\tau+1}}\Big)\widetilde H_k\Big(\frac{y}{\sqrt{\tau+1}}\Big)e^{-y^2}dy.
     \end{split}
\end{equation}
where $\widetilde H_n(x)$ is the normalized Hermite polynomial defined by $\widetilde H_n(x)=e^{x^2/2}\widehat H_n(x)$. In fact, for $A_s^{\pm}(\xi)$, $B_s^{\pm}(\xi)$ given in \eqref{DNNcoefficients}, we have
	\begin{equation}
		\begin{split}
			C_{nk}^{\pm} &=-\frac{(1\pm e^{-2})\lambda}{2(2\pi)^{\frac{3}{2}}(s+1)}\sum\limits_{j=0}^s\alpha^3_jI_{nk}\Big(\frac{2\pi^2}{\alpha^2_j\lambda^2}\Big),\quad \\
			M_{nk}^{\pm} &=-\sqrt{\frac{\pi}{2}}\frac{1\pm e^{-2}}{(s+1)\lambda}\sum\limits_{j=0}^s\alpha_jI_{nk}\Big(\frac{2\pi^2}{\alpha^2_j\lambda^2}\Big).
		\end{split}
	\end{equation}

	Next, we present formulas for the calculation of the integrals $I_{nk}(\tau)$. Given any scaling factor $\lambda$, scaled Hermite polynomial $\widetilde H_n(\lambda y)$ can be represented by $\widetilde H_n(y)$ as follows
	\begin{equation}\label{H_n(ty)_HZK}
		\widetilde H_n(\lambda y)=\sum_{k=0}^{n}h_{n,k}(\lambda)\widetilde H_k(y),
	\end{equation}
	where $\{h_{n,k}(\lambda)\}$ can be calculated via recurrence formulas \eqref{h_{k,m}_DTGS}. Therefore,
	\begin{equation}
		\begin{aligned}
			I_{nk}(\tau)&=\frac{1}{\sqrt{\tau+1}}\int_{-\infty}^{\infty}\widetilde H_n\Big(\frac{y}{\sqrt{\tau+1}}\Big)\widetilde H_k\Big(\frac{y}{\sqrt{\tau+1}}\Big)e^{-y^2}dy\\
			&=\frac{1}{\sqrt{\tau+1}}\sum_{i=0}^{n}\sum_{j=0}^{k}h_{n,i}\Big(\frac{1}{\sqrt{\tau+1}}\Big)h_{k,j}\Big(\frac{1}{\sqrt{\tau+1}}\Big)\int_{-\infty}^{\infty}\widetilde H_i(y)\widetilde H_j(y)e^{-y^2}dy\\
			&=\frac{1}{\sqrt{\tau+1}}\sum_{i=0}^{\min\{n,k\}}h_{n,i}\Big(\frac{1}{\sqrt{\tau+1}}\Big)h_{k,i}\Big(\frac{1}{\sqrt{\tau+1}}\Big).
		\end{aligned}
	\end{equation}
	
	Next, we derive recurrence formulas for the computation of the coefficients $\{h_{nk}(\lambda)\}$. We drop the explicit dependence on $\lambda$ without confusion in the following derivation. By the definition of $\widetilde H_n(y)$ and the recurrence formula \eqref{H_nDS}, we have
	\begin{equation}\label{H_nDS3}
		\sqrt{2(n+1)}\widetilde H_{n+1}(\lambda y)=2\lambda y\widetilde H_n(\lambda y)-\sqrt{2n}\widetilde H_{n-1}(\lambda y),\quad n\ge 1.
	\end{equation}
	Substituting the expansion \eqref{H_n(ty)_HZK} into \eqref{H_nDS3} gives for $n\geq 1$
	\begin{equation}\label{exprecurrence}
		\sqrt{2(n+1)}\sum_{k=0}^{n+1}h_{n+1,k}(\lambda)\widetilde H_k(y)=2\lambda y\sum_{k=0}^{n}h_{n,k}(\lambda)\widetilde H_k(y)-\sqrt{2n}\sum_{k=0}^{n-1}h_{n-1,k}(\lambda)\widetilde H_k(y).
	\end{equation}
	Noting that
	\begin{equation}
		\widetilde H_{1}(y)=\sqrt{2}y\widetilde H_0(y),\quad 2y\widetilde H_k(y)=\sqrt{2(k+1)}\widetilde H_{k+1}(y)+\sqrt{2k}\widetilde H_{k-1}(y)\ ,\quad k\ge 1,
	\end{equation}
	direct calculation from \eqref{exprecurrence} gives
	\begin{equation*}
		\begin{split}
			& 2\lambda y\sum_{k=0}^{n}h_{n,k}(\lambda)\widetilde H_k(y) \\
   =&\lambda\sum_{k=1}^{n}h_{n,k}(\lambda)\left[\sqrt{2(k+1)}\widetilde H_{k+1}(y)+\sqrt{2k}\widetilde H_{k-1}(y)\right]+2ayh_{n,0}(\lambda)\widetilde H_0(y)\\
			=&a\sum_{k=0}^{n}\sqrt{2(k+1)}h_{n,k}(\lambda)\widetilde H_{k+1}(y)+a\sum_{k=1}^{n}\sqrt{2k}h_{n,k}(\lambda)\widetilde H_{k-1}(y)\\
			=&a\sum_{k=1}^{n+1}\sqrt{2k}h_{n,k-1}(\lambda)\widetilde H_{k}(y)+a\sum_{k=0}^{n-1}\sqrt{2(k+1)}h_{n,k+1}(\lambda)\widetilde H_{k}(y).
		\end{split}
	\end{equation*}
	Therefore, \eqref{exprecurrence} can be rearranged into
	\begin{equation*}
		\begin{aligned}
			&\sqrt{2(n+1)}\sum_{k=0}^{n+1}h_{n+1,k}\widetilde H_k(y)\\
			=&\lambda\sum_{k=1}^{n+1}\sqrt{2k}h_{n,k-1}(\lambda)\widetilde H_{k}(y)+\lambda\sum_{k=0}^{n-1}\sqrt{2(k+1)}h_{n,k+1}(\lambda)\widetilde H_{k}(y) \\
    & -\sqrt{2n}\sum_{k=0}^{n-1}h_{n-1,k}(\lambda)\widetilde H_k(y)\\
			=&[\sqrt{2}\lambda h_{n,1}(\lambda)-\sqrt{2n}h_{n-1,0}(\lambda)]\widetilde H_{0}(y)+\lambda\sqrt{2n}h_{n,n-1}(\lambda)\widetilde H_n(y) \\
   & +\lambda\sqrt{2(n+1)}h_{n,n}(\lambda)\widetilde H_{n+1}(y)\\
			&+\sum_{k=1}^{n-1}[\lambda\sqrt{2k}h_{n,k-1}(\lambda)+\lambda\sqrt{2(k+1)}h_{n,k+1}(\lambda)-\sqrt{2n}h_{n-1,k}(\lambda)]\widetilde H_k(y).\\
		\end{aligned}
	\end{equation*}
	Matching the coefficients on both sides of the above equation gives us
	\begin{equation}\label{h_{k,m}_DTGS}
		\begin{aligned}
			h_{n+1,0}(\lambda)&=\sqrt{\frac{1}{n+1}}\lambda h_{n,1}(\lambda)-\sqrt{\frac{n}{n+1}}h_{n-1,0}(\lambda),\\
			h_{n+1,k}(\lambda)&=\lambda\sqrt{\frac{k+1}{n+1}}h_{n,k+1}(\lambda)-\sqrt{\frac{n}{n+1}}h_{n-1,k}(\lambda)+\lambda\sqrt{\frac{k}{n+1}}h_{n,k-1}(\lambda), \\
    \qquad & \mbox{for} \quad 1\le k \le n-1,\\
			h_{n+1,k}(\lambda)&=\lambda\sqrt{\frac{k}{n+1}}h_{n, k-1}(\lambda),\quad k=n,n+1,
		\end{aligned}
	\end{equation}
	for all $n\ge1$, while the initial values are given by
	\begin{equation}\label{h_{k,m}_DTCZ}
		h_{0,0}(\lambda)=1,\quad h_{1,0}(\lambda)=0,\quad h_{1,1}(\lambda)=\lambda.
	\end{equation}
	By induction,  $h_{n,k}(\lambda)$ has explicit formula for all $k=0,1,\cdots,n$
	\begin{equation}
		h_{n,k}(\lambda)=\begin{cases}
			0,  &  \ n-k=2s+1,\\
			\displaystyle \sqrt{\frac{n!}{2^{n-k}k!}}\frac{1}{s!}\lambda^k(\lambda^2-1)^s,  &  \ n-k=2s.
		\end{cases}
	\end{equation}

\end{appendix}

\end{document}